\documentclass[a4paper,fleqn]{cas-sc}

\usepackage[numbers]{natbib}
\usepackage{amsmath,amsfonts}
\usepackage{algorithmic}
\usepackage{algorithm}
\usepackage{array}
\usepackage{amssymb}
\usepackage[caption=false,font=normalsize,labelfont=sf,textfont=sf]{subfig}
\usepackage{textcomp}
\usepackage{stfloats}
\usepackage{url}
\usepackage{verbatim}
\usepackage{graphicx}
\usepackage{bm,bbm} 
\usepackage{booktabs} 
\usepackage{makecell} 
\usepackage{tabularx}
\usepackage{supertabular}
\usepackage{multirow}
\usepackage{bbding}
\usepackage{pifont}
\usepackage{array}
\usepackage[demo,export]{adjustbox}
\usepackage{colortbl}
\usepackage{hyperref} 
\newdefinition{rmk}{Remark}
\newproof{pf}{Proof}
\newproof{pot}{Proof of Theorem \ref{thm2}}
\def\tsc#1{\csdef{#1}{\textsc{\lowercase{#1}}\xspace}}
\tsc{WGM}
\tsc{QE}
\tsc{EP}
\tsc{PMS}
\tsc{BEC}
\tsc{DE}

\begin{document}
\let\WriteBookmarks\relax
\def\floatpagepagefraction{1}
\def\textpagefraction{.001}
\shorttitle{Multi-view fusion regularized clustering }
\shortauthors{X.R. Xing et~al.}

\title [mode = title]{Multi-View Clustering Meets Heterogenous Data: A Fusion Regularized Method}                      
\tnotemark[1]

\tnotetext[1]{This work was supported by the National Natural Science Foundation of China under Grant 12371322 and 12371306,  the Innovative Talent Project of Karamay under Grant XQZX20240078, and the Middle-aged and Young Teachers' Basic Ability Promotion Project of Guangxi (2025KY0025).}


\author[1]{Xiangru Xing}
\ead{21118025@bjtu.edu.cn}
\credit{Writing--original draft, Methodology, Data curation,  Visualization, Conceptualization}
\affiliation[1]{organization={School of Mathematics and Statistics, Beijing Jiaotong University},
                city={Beijing},
                postcode={100044}, 
                country={China}}

\author[2]{Yan Li}
\ead{liyan2010@uibe.edu.cn}
\credit{Writing--review,  Validation, Supervision}
\affiliation[2]{organization={School of Insurance and Economics, University of International Business and Economics},
	city={Beijing},
	postcode={100029}, 
	country={China}}
\author[1]{Xin Wang}
   
\ead{98930264@bjtu.edu.cn}
\cormark[1]
\credit{Writing--review \& editing, Methodology, Supervision, Resources}

\author[3]{Huangyue Chen}
\ead{hychen2024@gxu.edu.cn}
\credit{Writing– review \& editing, Methodology, Validation, Conceptualization}
\affiliation[3]{organization={School of Mathematics and Information Science, Guangxi University},
                city={Nanning},
                postcode={530004}, 
                country={China}}

\author[4]{Xianchao Xiu}
\ead{xcxiu@shu.edu.cn}
\credit{Writing--review \& editing, Methodology,  Formal analysis, Conceptualization}
\affiliation[4]{organization={School of Mechatronic Engineering and Automation, Shanghai University},
	city={Shanghai},
	postcode={200444},  
	country={China}}
	
\cortext[cor1]{Corresponding author}

\begin{abstract}
Multi-view clustering leverages consistent and  complementary information across multiple views to provide more comprehensive insights than single-view analysis. However, the heterogeneity and redundancy of multi-view data pose significant challenges to the existing clustering techniques. To  tackle these challenges effectively,
this paper proposes a novel multi-view fusion regularized clustering method with adaptive group sparsity, enabling discriminative clustering while  capturing informative features. Technically, for heterogeneous  multi-view data with mixed-type feature sets, different losses or divergence metrics are considered with a joint fusion penalty to obtain consistent cluster structures.
Moreover, the non-convex group sparsity consisting of inter-group sparsity and intra-group sparsity is utilized to eliminate redundant features, thereby enhancing the robustness. Furthermore, we develop an effective alternating direction method of multipliers (ADMM), where all subproblems can be solved in closed form. Extensive numerical experiments on  real data validate the superior performance of our presented method in clustering accuracy and feature selection.
\end{abstract}



\begin{keywords}
Multi-view clustering \sep Heterogenous data \sep Group sparsity \sep Feature selection \sep Alternating direction method of multipliers
\end{keywords}

\maketitle

\section{Introduction}
In the fields of image and audio processing, multi-view data are ubiquitous \cite{zhou2024survey}. For instance, images  can be characterized from multiple views such as the texture features extracted by local binary pattern (LBP), the shape features obtained through the histogram of oriented gradient (HOG), and the color features calculated via color moment (CMT).
These different representations often exhibit inherent relationships \cite{fang2023comprehensive}. A large number of studies have demonstrated that multi-view learning can fully exploit the diversity and complementarity  of these views for better efficiency and robustness \cite{yang2024methods}. Currently, it has become a very promising field, with wide applications in social recommendation \cite{ZHOU202412}, biometric recognition \cite{song2024study}, remote sensing \cite{Qi10950}, and information retrieval \cite{Sun2025113}.

Multi-view data collected from different sources usually have two prominent characteristics: \emph{heterogeneity} and \emph{redundancy}. For the former, it refers to data originating from diverse sources, dimensions, or observational perspectives that exhibit significant disparities in structure, scale, distribution, or semantic representation \cite{XUE2024,fogel2024integrated,HAO202412}. Additionally, its samples may contain multiple feature sets of varying types, such as continuous, count, binary, and categorical. This kind of heterogenous multi-view data is also termed \emph{mixed multi-view data} \cite{baker2020feature,Wang2021}. 
For example, in integrative genomic analyses, small RNA expression or gene  data derived from sequencing assay is typically zero-inflated  right-skewed continuous or count-valued, whereas DNA methylation is represented as proportion-valued data. These data types frequently deviate from Gaussian assumptions, rendering conventional Euclidean distance metrics suboptimal for analysis due to their sensitivity to non-Gaussian distributions. 
Crucially, while extensive methods (e.g., \cite{choi2019convex,chu2021adaptive,sun2025resistant}) have addressed such non-Gaussian distributions within individual data types, joint analytical approaches for mixed multi-view data remain underdeveloped. For the latter, one classical technique for reducing the dimensionality is to perform principal component analysis (PCA)  \cite{Ghosh2002},  $t$-distributed stochastic neighbor embedding ($t$-SNE) \cite{Van2008}, non-negative matrix factorization (NMF) \cite{XIANG2024} or projection \cite{fang2024joint} before clustering. However, such methods may lead to results that do not directly illustrate the importance of features. Moreover, in practical scenarios, determining the ranking of factors  is acknowledged to be extremely challenging and usually has a profound effect on the final outcomes. \cite{zhao2025multi}.

\begin{figure*}[t]
	\centering
	\includegraphics[width=0.92\textwidth]{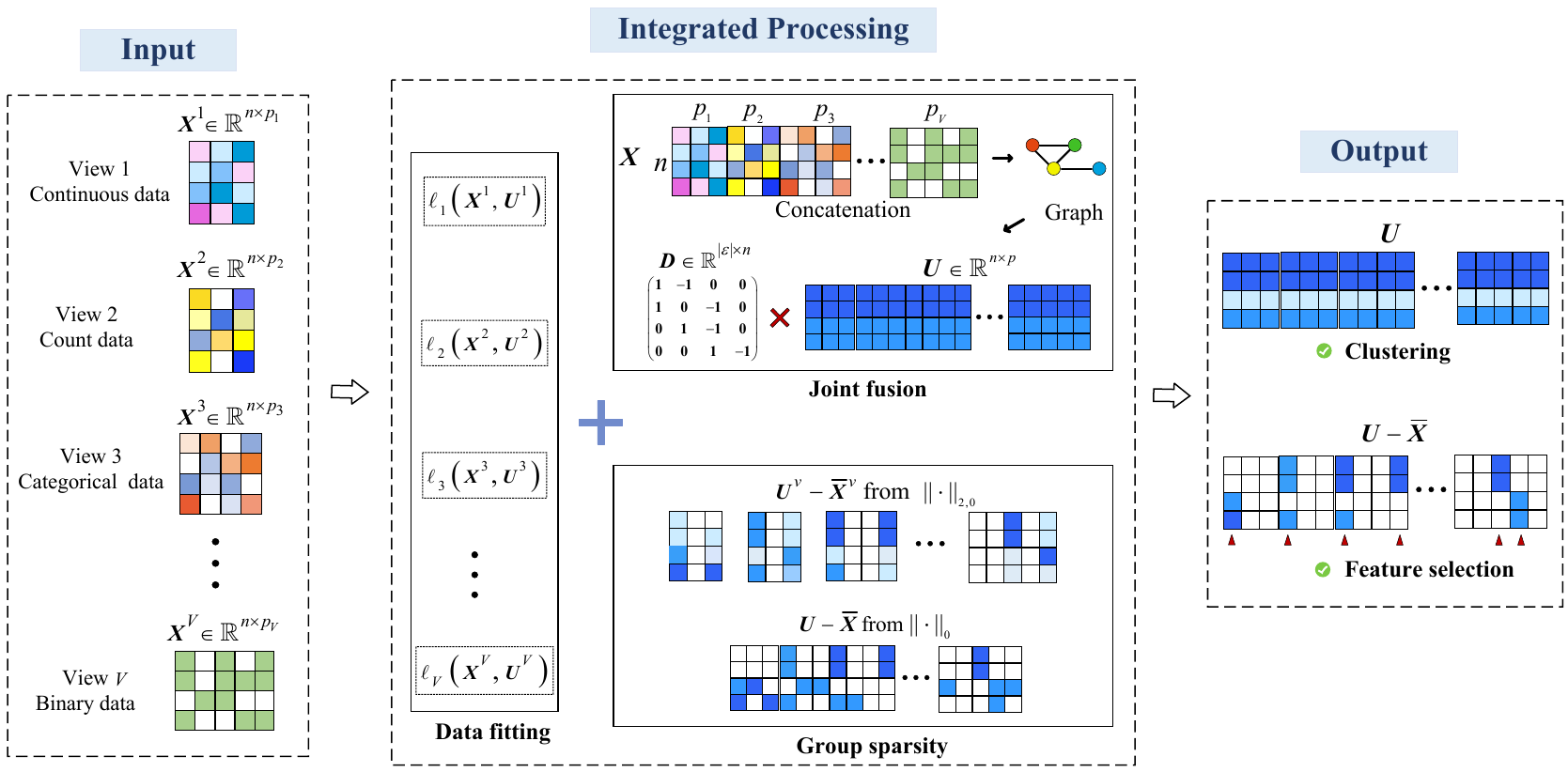}
	\caption{A framework of our proposed method. It integrates heterogenous multi-view data fitting, joint fusion, and feature selection into a unified framework to mine the global clustering structure and implement feature selection. Group sparsity captures both global and local feature structures, leading to different selection results compared to each individual sparsity term.}
	\label{fig1}
\end{figure*}

In this paper, we attempt to exploit heterogenous multi-view data to cluster common samples more effectively and select informative features to distinguish the inherent group structures.  Specifically, we present an innovative group-sparsity-based multi-view fusion regularized clustering method that can integrate mixed-type data via distinct data-specific losses, cluster co-samples through a joint fusion penalty, and select locally relevant features. The general framework for our method is shown in Figure \ref{fig1}.  Notably, in contrast to other prevalent multi-view clustering techniques, such as multi-view $K$-means clustering \cite{GUO202412}, multi-view subspace clustering  \cite{ZHOU2024}, and multi-view spectral clustering \cite{xie2024consistent},  our methodology  does not require prior knowledge regarding the number of clusters. To address non-Gaussian scenarios more effectively, a comprehensive loss function is formulated based on the characteristics of each view. Furthermore, inter-group sparsity with $l_{2,0}$-norm \cite{zhang2022structured,zhu2025joint}, a component of group sparsity, is devised to provide a more intuitive explanation than its approximate method called Lasso in \cite{Zou2005,Wang2021}. As another integral part of the group sparsity, the intra-group sparsity with $l_{0}$-norm \cite{chen2022nonconvex,tang2024sparse}, is ubiquitous in  real-world applications. Therefore, we develop a group sparsity term, i.e., $l_{2,0}$-norm plus $l_{0}$-norm, that aims to select relevant features with useful data points. 
However, the non-convexity and non-smoothness of group sparsity bring great challenges to model analysis and  numerical optimization.

Compared with the existing works, our main contributions are outlined below.
\begin{enumerate}[1)]
	\item For heterogenous  data, we present a general multi-view clustering model with inter-group sparsity and intra-group sparsity, which has not been considered before. This model can simultaneously and efficiently implement clustering and feature selection,  and is therefore more robust than traditional clustering models. 
	\item For the proposed non-convex and non-smooth model, we develop an ADMM-based  optimization algorithm. Fortunately, all subproblems admit closed-form solutions, thus enabling  efficient solving.
	\item Numerical experiments on real data demonstrate that our proposed method achieves superior performance compared to both state-of-the-art fusion regularized clustering approaches and other advanced clustering techniques.
\end{enumerate}

The overview of the remaining article  is shown below.  Section \ref{section2}  briefly reviews the related work. Section \ref{section3} provides our multi-view clustering model and optimization algorithm. The experimental results are detailed in  Section \ref{section5},  with final conclusions are presented in Section \ref{section6}.

\section{Preliminaries}\label{section2}
\subsection{Notations}

For a given matrix $\bm A\in \mathbb{R}^{n \times p}$, $\bm A_{i\cdot}$, $\bm A_{\cdot j}$, and $\bm A_{ij}$ denote its $i$-th row, $j$-th column, and $(i, j)$-th element, respectively. $\|\bm A\|_{F}=(\sum_{i=1}^{n}\sum_{j=1}^{p}\bm A^{2}_{ij})^{1/2}$ denotes the Frobenius norm. $\|\bm A\|_{2,1}=\sum_{i=1}^{n} \|\bm A_{i\cdot}\|_{2}$, $\|\bm A\|_{2,0}=\sum_{j=1}^{p} \mathbb{I}(\|\bm A_{\cdot j}\|_{2} \neq 0)$, $\|\bm A\|_{0}=\sum_{i=1}^{n}\sum_{j=1}^{p} \mathbb{I}(\bm A_{ij} \neq 0)$, where $\mathbb{I}(\cdot)$ is the indicator function. $\|\bm A_{i\cdot}\|_{q}=(\sum_{j=1}^{p}\bm A^{q}_{ij})^{1/q}$ denotes the $l_q$-norm of vector $\bm A_{i\cdot}$ with $q>0$. 
Besides,  $\bm e_n$ is an  $n$-dimensional column vector with components equal to  $1$. If not specified, $\bm I$ is  the identity matrix of appropriate dimension.

Let $f : \mathbb{R}^{n \times p} \to (-\infty, +\infty]$ be a proper closed function. For matrix $\bm A \in \mathbb{R}^{n \times p}$ and constant $\tau>0$, the proximal mapping of $f$ is formally defined by
\begin{eqnarray}
		\text{Prox}_{\tau f}(\bm A) = \underset{\bm U\in\mathbb{R}^{n \times p}}{\arg \min}\left\{\tau f(\bm U) + \frac{1}{2}\|\bm U - \bm A\|_F^2\right\}.\label{YY0}
\end{eqnarray} 

%
	%
		\begin{table}[width=.56\linewidth,cols=2,pos=t]
			\caption{Proximal mappings for selected functions.}
			\label{tab0}
			\begin{tabular*}{\tblwidth}{@{} LL@{} }
			\toprule
				$f(\cdot)$ &$\text{Prox}_{\tau f}(\bm A)$ \\ 
				\midrule
				$\|\cdot\|_{2,1}$ &$[\text{Prox}_{\tau\|\cdot\|_{2,1}}(\bm A)]_{i\cdot}=\left(1- \frac{\tau}{\max\{\|\bm A_{i\cdot}\|_2, \tau\}}\right) \mathbf{A}_{i\cdot}$   \\ 
				\midrule
				$\|\cdot\|_{2,0}$ &$[\text{Prox}_{\tau\|\cdot\|_{2,0}}(\bm A)]_{\cdot j}=\begin{cases}
					\bm A_{\cdot j}, & \|\bm A_{\cdot j}\|_2 > \sqrt{2\tau} \\
					\bm A_{\cdot j} \text{ or } 0, & \|\bm A_{\cdot j}\|_2 = \sqrt{2\tau} \\
					0, & \|\bm A_{\cdot j}\|_2 < \sqrt{2\tau}
				\end{cases}$  \\ 
				\midrule
				$\|\cdot\|_{0}$ &$[\text{Prox}_{\tau\|\cdot\|_{0}}(\bm A)]_{ij}=\begin{cases}
					\bm A_{ij}, & |\bm A_{ij}| > \sqrt{2\tau} \\
					\bm A_{ij} \text{ or } 0, & |\bm A_{ij}| = \sqrt{2\tau} \\
					0, & |\bm A_{ij}| < \sqrt{2\tau}
				\end{cases}$  \\ 
				\bottomrule
			\end{tabular*}
		\end{table}

	We list the proximal mappings involved in this paper in Table \ref{tab0}. More functions and associated proximal mappings can be found in \cite{beck2017first}.
		
	\subsection{Related Work}
	The fusion regularized clustering (alternatively termed sum-of-norms clustering, convex clustering, or clusterpath) was initially proposed by
	\cite{Pelckmans2005},  whose fundamental model is 
	\begin{eqnarray}\label{YY1}
			\min_{\bm U}\;\frac{1}{2}\sum_{i=1}^n\|\bm X_{i\cdot}- \bm U_{i\cdot}\|_2^2+\eta\sum_{i<i^\prime}\omega_{\iota}\|\bm U_{i\cdot}-\bm U_{i^\prime\cdot}\|_q,
	\end{eqnarray}
	where $\bm X\in\mathbb{R}^{n\times p}$ denotes the data matrix containing $n$ observations with $p$ features, $\bm U\in\mathbb{R}^{n\times p}$ is an centroid matrix, $\eta$ is a positive tuning parameter that balances  model fit and  cluster quantity. By adjusting $\eta$, the number of clusters can be automatically derived  without  prior specification \cite{Lindsten2011}. In addition, $\|\cdot\|_q$ is the $l_q$-norm  with $q\in\{1, 2,\infty\}$ commonly considered. For example,  $l_1$-norm loss is utilized in the presence of outliers. Given a positive constant $\phi$, the fusion weight $\omega_{\iota}$ is defined as
	\begin{eqnarray}
				\omega_{\iota}=\begin{cases}					\;\exp\big(-\phi \|\bm X_{i\cdot}-\bm X_{i^\prime\cdot}\|_2^2\big), &{\text{if}} \ (i, i^\prime)\in\epsilon,\\ 
					0, &{\text{otherwise}},
									\end{cases}
	\end{eqnarray} 
	where  $\epsilon=\cup_{\iota_1=1}^n\{ (i,i^\prime) \mid \bm X_{i^\prime\cdot}$ is among $\bm X_{i\cdot} \text{'s $\mathcal K$-nearest}$ neighbors, $1\le i < i^\prime \le n\}$.

	Although fusion regularized clustering  has desirable theoretical results and good computational performance \cite{Sun2021}, its clustering performance may degrade in high-dimensional cases with more features than observations. 
	To address this,  Pelckmans et al. \cite{Wang2018} introduced sparse convex clustering  incorporating an inter-group Lasso penalty, i.e.,
		\begin{eqnarray}
			\begin{aligned}\label{YY2}
				\min_{\bm U}\;\frac{1}{2}\sum_{i=1}^n\|\bm X_{i\cdot}- \bm U_{i\cdot}\|_2^2+\eta\sum_{i<i^\prime}\omega_{\iota}\|\bm U_{i\cdot}-\bm U_{i^\prime\cdot}\|_q+\beta\sum_{j=1}^p\bar{\beta}_j\|\bm U_{\cdot j}\|_2,
			\end{aligned}
		\end{eqnarray} 
		where tuning parameter $\beta>0$ regulates the sparsity level of  features and weight coefficient $\bar{\beta}_j$ determines the importance ratio of each feature. In this way, the cluster centers are shrunk towards zero and relevant features are selected. Building on this,  Chen et al.
		\cite{Chen2020} considered the intra-group sparse penalty with $l_1$-norm to further enhance the ability of feature selection. However, these models are only suitable for data with a mean of zero. Otherwise, the Lasso penalty terms will select incorrect features, leading to biased cluster centers.
		
		Moreover,    the above fusion regularized clustering models  often perform poorly with non-Gaussian data or lack robustness to noise \cite{feng2023review}. 
		To address these limitations, researchers have adopted alternative losses, such as Bernoulli log-likelihood \cite{choi2019convex},  Huber loss \cite{sun2025resistant}, and others \cite{chu2021adaptive,ma2023improved}, with the aim of improving the accuracy of data fitting. Recently,  Wang and Allen \cite{Wang2021} proposed an adaptive inter-group Lasso penalty for heterogeneous multi-view data with mixed feature types. This approach performs feature selection by shrinking  features toward their  loss-speciﬁc centers, which is given by 
		\begin{eqnarray}
			\begin{aligned}\label{YY3}
				\min_{\bm U^1, \bm U^2,\cdots, \bm U^V}\quad\sum_{v=1}^{V}\zeta_v \ell_v(\bm X^v, \bm U^v)+\eta\sum_{i<i^\prime}\omega_{\iota}\times \sqrt{\sum_{v=1}^{V}\|\bm U^v_{i\cdot}-\bm U^v_{i^\prime\cdot}\|_2^2}+\beta \sum_{v=1}^{V}\sum_{j=1}^{p_v}\bar{\beta}^v_j\|\bm U^v_{\cdot j}-\bar{\bm x}^v_j\cdot\bm e_n\|_2,
			\end{aligned}
		\end{eqnarray} 
		where $V$ represents the total number of views. For the $v$-th view, $\bm X^v\in\mathbb{R}^{n\times p_v}$ is the $v$-th data matrix with $n$ observations and $p_v$ features,  $\bm U^v\in\mathbb{R}^{n\times p_v}$ is the corresponding centroid, $\zeta_v$ is a user-specified weight with $\sum_{v=1}^{V}\zeta_v=1$,  the $v$-th loss function
		$\ell_v(\bm X^v, \bm U^v)$ is appropriately chosen  depending on the type of $\bm X^v$. In addition,  $\bar{\bm x}^v=\text{arg min}_{\bm U^v}~\ell_v(\bm X^v, \bm U^v)$ is a loss-specific center, with $\bar{\bm x}^v_j$ being its   $j$-th component. For instance, if $\ell_v(\bm X^v, \bm U^v)=\sum_{i=1}^{n}\|\bm X^v_{i\cdot}-\bm U^v_{i\cdot}\|_1$ for count data, then $\bar{\bm x}^v$ is the median of observations. See \cite[Table 1]{Wang2021}  for additional examples.
		It is worth noting that \eqref{YY3} is the first to address fusion regularized clustering for mixed multi-view data, capturing the inherent structure of diverse data types while offering greater generality.

		\section{Proposed Method}\label{section3}
		
		This section  starts with the presentation of our new model, followed by its optimization algorithm.
		
		\subsection{New Model}
		
		To deal with inherent  \emph{heterogeneity} and \emph{redundancy} in  multi-view data, we construct a multi-view fusion regularized clustering model with group sparsity  defined as
		\begin{eqnarray}
			\begin{aligned}\label{XX01}
				\min_{\bm U}\;\sum_{v=1}^{V}\zeta_v \ell_v(\bm X^v, \bm U^v)+\eta\sum_{\iota\in\epsilon}\omega_\iota\|(\bm D\bm U)_{\iota\cdot}\|_2+
\beta \left[(1-\theta)\sum_{v=1}^{V}\|(\bm U^v-\bar{\bm X}^v)\bar{\bm \beta}^v\|_{2,0}+\theta\|\bm U-\bar{\bm X}\|_0\right].
			\end{aligned}
		\end{eqnarray}
		where $\bm U=[\bm U^1,\;\bm U^2,\;\ldots,\;\bm U^V]\in\mathbb{R}^{n\times p}$ is the whole centroid matrix with $\bm U^v\in\mathbb{R}^{n\times p_v}$. 
		Mathematically, the first two terms are the matrix representations of the corresponding terms in \eqref{YY3}. Next, we will explain each item in detail.
		\begin{itemize}
			\item The first term plays the role of data fitting. Loss function $\ell_v(\bm X^{v}, \bm U^{v})$ is convex and possibly non-smooth for the $v$-th data view whose type is chosen appropriately according to the characteristics of  data $\bm X^v$. Here, we suppose that the non-smooth loss $\ell_v(\bm X^{v}, \bm U^{v})$ is distance-based such as Manhattan distance, Chebychev distance, and square-root-loss.
			
			\item The second term facilitates the rows' differences of concatenated centroid $[\mathbf{U}^{1},\; \mathbf{U}^{2},\;\ldots,\; \mathbf{U}^{V}]$ to shrink towards zero, thereby inducing a clustering behavior. In detail,  $(\bm D\bm U)_{\iota\cdot}=\bm U_{i\cdot}-\bm U_{i^\prime\cdot}$ for edge $\iota=(i, i^\prime)\in\epsilon$ and the directed difference matrix $\bm D\in\mathbb{R}^{|\epsilon|\times n}$. We use $l_2$-norm to promote the fusion of the whole rows of similar observations simultaneously and to ensure rotation invariance.  Under the action of the sparse fusion weight $\omega_{\iota}$,
			this term enforces the group structure of the $i$-th row of $\mathbf{U}^{v}$ to be  identical  for  $\forall v\in[V]$.  We consider that $\bm X_{i\cdot}$ and $\bm X_{i^\prime\cdot}$ belong to the same cluster if $\mathbf{U}_{i\cdot}=\mathbf{U}_{i^\prime\cdot}$.  
			Unlike fusion weight $\omega_{\iota}$ in convex clustering  models \eqref{YY1} and \eqref{YY2}, we introduce 
		\begin{eqnarray}
				\omega_{\iota}=\begin{cases}					\;\exp\Big(-\phi \sum_{v=1}^V\sum_{j=1}^{p_v}\frac{d_{ii^\prime j}^{v}}{\sum_{k=1}^V p_v}\Big), &{\text{if}} \ (i, i^\prime)\in\epsilon,\\ 
					0, &{\text{otherwise}},
				\end{cases}
		\end{eqnarray}
		\vspace{-0.1em}
		based on the Gower distance	 for clustering heterogenous multi-view data \cite{wangchamhan2017efficient,akay2018clustering}, where 
		\begin{eqnarray}
				d_{i i^\prime j}^{v}=\frac{|\mathbf{X}_{i j}^{v}-\mathbf{X}_{i^\prime j}^{v}|}{\max _{i,  i^\prime}|\mathbf{X}_{i j}^{v}-\mathbf{X}_{ i^\prime j}^{v}|}.
		\end{eqnarray}
			
			\item The last two terms are enforced to achieve feature selection. Here, $\bar{\bm X}^v=[\bar{\bm x}^v_{ 1}\cdot\bm e_n,\; \bar{\bm x}^v_{ 2}\cdot\bm e_n,\; \cdots, \;\bar{\bm x}^v_{ p_{v}}\cdot\bm e_n]$ with $\bar{\bm x}^v$  is consistent with the definition in \eqref{YY3} representing a loss-specific center, and $\bar{\bm X}$ is the horizontal concatenation matrix of  $\{\bar{\bm X}^v\}_{v\in[V]}$. By subtracting   $\bar{\bm X}^v$,  we perform centering on the centroid  $\bm U^v$, which shifts the feature means to zero and removes any bias.  The parameter $\theta\in[0, 1]$  balances the inter-group sparsity and intra-group sparsity,  while $\bar{\bm \beta}^v=\textup{Diag}(\bar{\beta}^v_1, \bar{\beta}^v_2, \bar{\beta}^v_3, \cdots, \bar{\beta}^v_{p_v})$ (user-specified, $\bar{\beta}^v_j>0$ for $j\in[p_v]$)  adaptively determines the importance of each feature.
			The impact of tuning parameters $\eta$, $\beta$, and $\theta$ will be explained in Section \ref{section5}.
		\end{itemize}
		
		We would also like to point out that although the Lasso method can reveal the sparse structure of solutions, it may lead to biased estimates as illustrated in \cite{Fan2001}. Different from the multi-view fusion regularized clustering with Lasso penalty discussed in \cite{Chen2020,Wang2021}, our proposed model considers the group sparsity that is the combination of $l_{2,0}$-norm and $l_0$-norm. In feature selection, the group sparsity provides a more intuitive explanation because $l_0$-norm is the primordial method for characterizing sparsity. In addition, group sparsity has great flexibility by adjusting the parameters $\beta$ and $\theta$.

		\subsection{Optimization Algorithm}
		
		Due to its non-differentiability with respect to the variable $\bm U$, problem \eqref{XX01} poses   great  difficulties to solution. 
		Therefore, it is necessary to develop an effective algorithm.

		We  introduce four slack variables $\bm B^{v}\in\mathbb{R}^{n\times p_v}$, $\bm E\in\mathbb{R}^{n\times p}$, $\bm F^{v}\in\mathbb{R}^{n\times p_v}$, and $\bar{\bm F}^{v}\in\mathbb{R}^{n\times p_v}$  for ease of solving, and the following notations   
	\begin{eqnarray}
		\begin{aligned}
		&	L(\bm U^{v})= \zeta_{v} \ell_{v}\left(\bm X^{v}, \bm U^{v}\right),  \\
		&	q(\bm B^{v})=\zeta_v\ell_v(\bm B^{v}),~p(\bm E)=\eta\sum_{\iota\in\epsilon}w_\iota\|(\bm E)_{\iota\cdot}\|_2, 	
			\\
		&	r_1(\bm F^{v})=\beta(1-\theta)\|\bm F^v\bar{\bm \beta}^v\|_{2,0}, ~r_2(\bar{\bm F^{v}})=\beta\theta\|\bar{\bm F}^v\|_0.
		\end{aligned}
\end{eqnarray}

Since the loss function  $\ell_v(\bm X^{v},\bm U^{v})$  may be either non-differentiable or differentiable, for $\forall v\in[V]$. For convenience, we will discuss it in two cases.

(Case I)
When $\ell_{v}\left(\bm X^{v}, \bm U^{v}\right)$ is not differentiable, problem \eqref{XX01} for each view can be equivalently transformed to
		\begin{eqnarray}
			\begin{aligned}\label{XX03}
			\min_{\bm U^{v}, \bm B^{v},\bm E, \bm F^{v},  \bar{\bm F}^{v}}~&\;  L(\bm U^{v})+q(\bm B^{v})+p(\bm E)+r_1(\bm F^{v})
			+r_2(\bar{\bm F}^{v})
			\\
			\textrm{s.t.}~~~~~~~~~&\; \bm{ DU}=\bm E,  \quad  \bm X^{v}- \bm U^{v}=\bm B^{v},
			\\
			&\; \bm U^{v}-\bar{\bm X}^{v}=\bm F^{v},\quad \bm U^{v}-\bar{\bm X}^{v}=\bar{\bm F}^{v},
			\end{aligned}
		\end{eqnarray}
where function $p(\cdot)$  is non-separable with respect to each view $v$.

		Let $\sigma> 0$, the augmented Lagrangian function of \eqref{XX03} is given by
		\begin{eqnarray}
			\begin{aligned}\label{XX06}
				&\mathcal{L}_\sigma(\bm U^{v}, \bm B^{v}, \bm E, \bm F^{v},  \bar{\bm F}^{v};  \bm Q^{v}, \bm P, \bm G^{v}, \bar{\bm G}^{v})
				\\
				&=L(\bm U^{v})+q(\bm B^{v})+p(\bm E)+r_1(\bm F^{v})+r_2(\bar{\bm F}^{v})
				+\langle\bm Q^{v}, \bm U^{v}+\bm B^{v}-\bm X^{v}\rangle
				+\frac{\sigma}{2}\|\bm U^{v}+\bm B^{v}-\bm X^{v}\|_F^2
				\\
				&+\langle\bm P, \bm{ DU}-\bm E\rangle+\frac{\sigma}{2}\|\bm{ DU}-\bm E\|_F^2,
				+\langle\bm G^{v}, \bm U^{v}-\bar{\bm X}^{v}-\bm F^{v}\rangle+\frac{\sigma}{2}\|\bm U^{v}-\bar{\bm X}^{v}-\bm F^{v}\|_F^2
				\\
				&+\langle\bar{\bm G}^{v}, \bm U^{v}-\bar{\bm X}^{v}-\bar{\bm F}^{v}\rangle+\frac{\sigma}{2}\|\bm U^{v}-\bar{\bm X}^{v}-\bar{\bm F}^{v}\|_F^2,
			\end{aligned}
		\end{eqnarray}
		where $\bm Q^{v}\in\mathbb{R}^{n\times p_v}$, $\bm P\in\mathbb{R}^{n\times p}$, $\bm G^{v}\in\mathbb{R}^{n\times p_v},
		 \;\text{and} \;\bar{\bm G}^{v}\in\mathbb{R}^{n\times p_v}$ are the Lagrange multipliers.

(Case II) When $\ell_{v}\left(\bm X^{v}, \bm U^{v}\right)$ is  differentiable, the corresponding view-specific subproblem in \eqref{XX01} can be equivalently reformulated as
\begin{eqnarray}
	\begin{aligned}\label{XX33}
		\min_{\bm U^{v}, \bm E, \bm F^{v},  \bar{\bm F}^{v}}~&\;  L(\bm U^{v})+p(\bm E)+r_1(\bm F^{v})+r_2(\bar{\bm F}^{v})
		\\
		\textrm{s.t.}~~~~~~&\; \bm{ DU}=\bm E, 
		\bm U^{v}-\bar{\bm X}^{v}=\bm F^{v}, \bm U^{v}-\bar{\bm X}^{v}=\bar{\bm F}^{v},
	\end{aligned}
\end{eqnarray}
whose augmented Lagrangian function  is
\begin{eqnarray}
	\begin{aligned}\label{XX66}
		&\mathcal{L}_\sigma(\bm U^{v},  \bm E, \bm F^{v},  \bar{\bm F}^{v};  \bm P, \bm G^{v}, \bar{\bm G}^{v})
		\\
		&=L(\bm U^{v})+p(\bm E)+r_1(\bm F^{v})+r_2(\bar{\bm F}^{v})+\langle\bm P, \bm{ DU}-\bm E\rangle+\frac{\sigma}{2}\|\bm{ DU}-\bm E\|_F^2,
		\\
		&+\langle\bm G^{v}, \bm U^{v}-\bar{\bm X}^{v}-\bm F^{v}\rangle+\frac{\sigma}{2}\|\bm U^{v}-\bar{\bm X}^{v}-\bm F^{v}\|_F^2,
		+\langle\bar{\bm G}^{v}, \bm U^{v}-\bar{\bm X}^{v}-\bar{\bm F}^{v}\rangle+\frac{\sigma}{2}\|\bm U^{v}-\bar{\bm X}^{v}-\bar{\bm F}^{v}\|_F^2.
	\end{aligned}
\end{eqnarray}

	As illustrated in \cite{chen2016direct}, the convergence of ADMM with multiple blocks cannot be guaranteed. Fortunately, \eqref{XX01} admits a parallel implementation via 2-block ADMM through the following variable partitioning
	\begin{eqnarray}
		\begin{aligned}
			\bm\Upsilon_1&=(\bm B^{v}, \bm F^{v},  \bar{\bm F}^{v}, \bm E),\\
			\bm\Upsilon_2&=(\bm Q^{v}, \bm G^{v}, \bar{\bm G}^{v}, \bm P).	
		\end{aligned}
	\end{eqnarray}
		To this end, we can alternately minimize the primal variables among $\bm\Upsilon_1$ and  $\bm U^{v}$, and then
		update the multipliers $\bm\Upsilon_2$. The order of $\bm B^{v}$, $\bm F^{v}$,  $\bar{\bm F}^{v}$, and  $\bm E$ in the first block $\bm\Upsilon_1$ is arbitrary, which has no impact on the accuracy of ADMM. 
		
		Given the current iteration point $(\bm\Upsilon_1^{(t)}, \bm U^{v(t)}, \bm\Upsilon_2^{(t)})$ at the $t$-th iteration (where $t$ denotes the iteration counter), then the next iteration point $(\bm\Upsilon_1^{(t+1)}, \bm U^{v(t+1)}, \bm\Upsilon_2^{(t+1)})$ is generated below.

		\subsubsection{Update $\bm B^{v(t+1)}$}
		
By incorporating  (\protect\ref{XX06}), $\bm B^{v}$-related problem (\protect\ref{XX03}) can be simplified to
		\begin{eqnarray}
			\begin{aligned}
				\min_{\bm B^{v}}\quad q(\bm B^{v})+\frac{\sigma}{2}\|\bm U^{v(t)}+\bm B^{v}-\bm X^{v}+\frac{1}{\sigma}\bm Q^{v(t)}\|_F^2.
			\end{aligned}
		\end{eqnarray}
		We can obtain that its closed-form solution  is 
		\begin{eqnarray}
			\begin{aligned}\label{XX07}
				\bm B^{v(t+1)}\in \textrm{Prox}_{\zeta_v\ell_v / \sigma}(\bm X^{v}-\bm U^{v(t)}+\frac{1}{\sigma}  \bm Q^{v(t)}), 
			\end{aligned}
		\end{eqnarray}
		where \textrm{Prox} is defined in (\ref{YY0}) and Table \ref{tab0}.
		
		\subsubsection{Update $\bm F^{v(t+1)}$}
		Likewise, the optimization problems (\protect\ref{XX03})-(\protect\ref{XX33}) concerning $\bm F^{v}$ can be reduced to
		\begin{eqnarray}
				\min_{\bm F^{v}}\quad r_1(\bm F^{v})+\frac{\sigma}{2}\|\bm U^{v(t)}-\bar{\bm X}^{v}-\bm F^{v}+\frac{1}{\sigma}\bm G^{v(t)}\|_F^2.
		\end{eqnarray}
		Its closed-form solution can be described by
		\begin{eqnarray}\label{XX08}
				\bm F^{v(t+1)}\in \textrm{Prox}_{ r_1/ \sigma}(\bm U^{v(t)}-\bar{\bm X}^v+\frac{1}{\sigma}\bm G^{v(t)}). 
		\end{eqnarray}

		\subsubsection{Update $\bar{\bm F}^{v(t+1)}$}
		The $\bar{\bm F}^{v}$-subproblem can be formulated using  (\protect\ref{XX06})-(\protect\ref{XX66}) as
		\begin{eqnarray}
				\min_{\bar{\bm F}^{v}}\quad r_2(\bar{\bm F}^{v})+\frac{\sigma}{2}\|\bm U^{v(t)}-\bar{\bm X}^{v}-\bar{\bm F}^{v}+\frac{1}{\sigma}\bar{\bm G}^{v(t)}\|_F^2,
		\end{eqnarray}
	whose closed-form solution  can be expressed as
		\begin{eqnarray}\label{XX09}
				\bar{\bm F}^{v(t+1)}\in \textrm{Prox}_{ r_2/ \sigma}(\bm U^{v(t)}-\bar{\bm X}^v+\frac{1}{\sigma}\bar{\bm G}^{v(t)}).
		\end{eqnarray}

		\subsubsection{Update $\bm E^{(t+1)}$}
		The $\bm E$-subproblem can be updated via
		\begin{eqnarray}
				\min_{\bm E}\quad p(\bm E)+\frac{\sigma}{2}\|\bm{ DU}^{(t)}-\bm E+\frac{1}{\sigma}\bm P^{(t)}\|_F^2,
		\end{eqnarray}
		which has a closed-form solution 
		\begin{eqnarray}\label{XX10}
				\bm E^{(t+1)}\in \textrm{Prox}_{ p/ \sigma}(\bm {DU}^{(t)}+\frac{1}{\sigma}\bm P^{(t)}).
		\end{eqnarray}

		\begin{algorithm}[t]
			\caption{ADMM for solving \eqref{XX01}}\label{alg:1}
			\textbf{Initialization: } Choose $(\bm\Upsilon_1^{(0)}, \bm U^{(0)}, \bm\Upsilon_2^{(0)})$ and $\sigma>0$ \\
			\textbf{while} not converged	\textbf{do}
			
			\hspace{0.5cm}\textbf{for} $t=0, 1, 2, \cdots$, \textbf{do}
			
			\hspace{0.5cm}\textbf{Step 1.}  Compute $\bm\Upsilon_1^{(t+1)}$:
			
			\hspace{1cm}	Update $\bm B^{v(t+1)}$ by (\ref{XX07})
			
			\hspace{1cm}    Update $\bm F^{v(t+1)}$ by (\ref{XX08})
			
			\hspace{1cm}    Update $\bar{\bm F}^{v(t+1)}$ by (\ref{XX09})
			
			\hspace{1cm}	Update $\bm E^{(t+1)}$ by (\ref{XX10})
			
			\hspace{0.5cm}\textbf{Step 2.}  Compute $\bm U^{(t+1)}$:
			
			\hspace{1cm}\textbf{if} $\ell_{v}$ is not differentiable, 
			
			\hspace{1cm} update
			$\bm U^{v(t+1)}$ by (\ref{XX12})
			
			\hspace{1cm}\textbf{else if} 
			
			\hspace{1cm} update
			$\bm U^{v(t+1)}$ by (\ref{XX15})
			
			\hspace{1cm}\textbf{end if}
			
			\hspace{0.5cm}\textbf{Step 3.}  Compute $\bm\Upsilon_2^{(t+1)}$:
			
			\hspace{1cm}Update $\bm Q^{v(t+1)}, \bm G^{v(t+1)}, \bar{\bm G}^{v(t+1)}, \bm P^{(t+1)}$ by (\ref{XX16})
			
			\hspace{0.5cm}\textbf{end for}\\
			\textbf{end while}\\
			\textbf{Output: } $(\bm\Upsilon_1^{(t+1)}, \bm U^{(t+1)}, \bm\Upsilon_2^{(t+1)})$
		\end{algorithm}
		
		\subsubsection{Update $\textbf U^{v(t+1)}$}
		(Case I) 	When $\ell_{v}\left(\bm X^{v}, \bm U^{v}\right)$ is not differentiable, 
		we can  obtain  the closed-form solution from (\ref{XX03}) and (\ref{XX06}) by solving
		\begin{eqnarray}
			\begin{aligned}\label{XX12}
			\bm U^{v(t+1)}&=\underset {\bm U^{v}} {\text{arg min}}\big\{\ell_v(\bm X^{v},\bm U^{v})
			+\frac{\sigma}{2}\|\bm U^{v}+\bm B^{v(t+1)}-\bm X^{v}+\frac{1}{\sigma}\bm Q^{v(t)}\|_F^2
			+\frac{\sigma}{2}\|\bm{ DU}-\bm E^{(t+1)}+\frac{1}{\sigma}\bm P^{(t)}\|_F^2
			\\
			&+\frac{\sigma}{2}\|\bm U^{v}-\bar{\bm X}^{v}-\bm F^{v(t+1)}+\frac{1}{\sigma}\bm G^{v(t)}\|_F^2
			+\frac{\sigma}{2}\|\bm U^{v}-\bar{\bm X}^{v}-\bar{\bm F}^{v(t+1)}+\frac{1}{\sigma}\bar{\bm G}^{v(t)}\|_F^2.
			\end{aligned}
		\end{eqnarray}

		(Case II) When $\ell_{v}\left(\bm X^{v}, \bm U^{v}\right)$ is  differentiable, we can acquire the closed-form solution from (\ref{XX33}) and (\ref{XX66})  as
		\begin{eqnarray}
			\begin{aligned}\label{XX15}
				\bm U^{v(t+1)}&=\big(3\bm I+\bm D^\top\bm D\big)^{-1}\big(\bm X^{v}-\bm B^{v(t+1)}
				-\frac{1}{\sigma}\bm Q^{v(t)}+\bm D^\top\bm E^{(t+1)}-\frac{1}{\sigma}\bm D^\top\bm P^{(t)}+2\bar{\bm X}^{v}
				+\bm F^{v(t+1)}\\
				&-\frac{1}{\sigma}\bm G^{v(t)}+\bar{\bm F}^{v(t+1)}
			-\frac{1}{\sigma}\bar{\bm G}^{v(t)}\big).
			\end{aligned}
		\end{eqnarray}

	\subsubsection{Update Multipliers}
	The multipliers in (\ref{XX06}) and (\ref{XX66}) are updated  as follows
		\begin{eqnarray}
		\begin{aligned}
		\label{XX16}
		\bm Q^{v(t+1)}&= \bm Q^{v(t)}+\sigma(\bm X^{v}-\bm U^{v(t+1)}-\bm B^{v(t+1)}),
\\
		\bm G^{v(t+1)}&= \bm G^{v(t)}+\sigma(\bm U^{v(t+1)}-\bar{\bm X}^v-\bm F^{v(t+1)}),
			\\
		\bar{\bm G}^{v(t+1)}&= \bar{\bm G}^{v(t)}+\sigma(\bm U^{v(t+1)}-\bar{\bm X}^v-\bar{\bm F}^{v(t+1)}),	
			\\
		\bm P^{(t+1)}&= \bm P^{(t)}+\sigma(\bm D\bm U^{(t+1)}-\bm E^{(t+1)}).
			\end{aligned}
	\end{eqnarray}

	The horizontally concatenate matrices (e.g.  $\bm F$, $\bar{\bm F}$, and $\bm U$) are divided into  matrices  of each view ($\bm F^{v}$, $\bar{\bm F}^{v}$, and $\bm U^{v}$)  for separate solution, so as to reduce some computational complexity without affecting the solution accuracy.  The specific iterative steps are shown in Algorithm \ref{alg:1}.

\subsection{Complexity Analysis}
The complexity  primarily stems from the update of  original variables  and  multipliers. As previously defined, $n$, $p_v$, $p$,  $|\epsilon|$, and $\mathcal{K}$ denote the number of observations, features in the $v$-th view, total features across  all views,   edges connecting observations, and nearest neighbors per  observation, respectively. The specific choice of $\mathcal{K}\in\{3,4,5,6,7,8,9,10\}$ in numerical experiments leads to $|\epsilon|>n$.

\begin{itemize}
	\item Time Complexity:
	The update process of the original variables focuses on   $\bm B^{v}$, $\bm F^{v}$, $\bar{\bm F}^{v}$, $\bm E$, and $\bm U^{v}$,  with corresponding  complexities of $\mathcal{O}(np_v)$, $\mathcal{O}(np_v)$, $\mathcal{O}(np_v)$, $\mathcal{O}(|\epsilon|p)$, and $\mathcal{O}(|\epsilon|p_v+n^3+n^2p_v)$. Here, for updating  $\bm B^{v}$, $\bm F^{v}$, and $\bar{\bm F}^{v}$, the time cost is dominated by matrices addition and subtraction. For updating $\bm E$, the  cost  depends on  computing  $\epsilon$ vectors subtraction $\bm {DU}=\big[\bm D_{ii^\prime}(\bm U_{i\cdot}-\bm U_{i^\prime\cdot})\big]_{(i,i^\prime)\in\epsilon}\in\mathbb{R}^{|\epsilon|\times p}$. Because the $\bm E$-subproblem cannot be separated across views,  updating  $\bm E$ and its corresponding multiplier $\bm P$ involves storing all features. The cost of updating $\bm U^{v}$ is highly dependent on computing matrix inverse $(3\bm I+\bm D^\top\bm D)^{-1}$ and matrix multiplication. Similarly, for the multipliers $\bm Q^{v}$, $\bm G^{v}$, $\bar{\bm G}^{v}$, and  $\bm P$, the costs are $\mathcal{O}(np_v)$, $\mathcal{O}(np_v)$, $\mathcal{O}(np_v)$, and $\mathcal{O}(|\epsilon|p)$. Overall, the  time complexity is $\mathcal{O}\big(t_1(|\epsilon|p+n^3+n^2p)\big)$, where $t_1$  represents the iteration count of Algorithm \ref{alg:1}.
	
	\item Space Complexity: The dominant space cost comes from storing matrices $\bm {E}\in\mathbb{R}^{|\epsilon|\times p}$, $\bm {D}\in\mathbb{R}^{|\epsilon|\times n}$,  and $\bm {P}\in\mathbb{R}^{|\epsilon|\times p}$ during the updates of  $\bm E$,  $\bm U^{k}$, and $\bm P$. Thus, the space complexity is $\mathcal{O}(|\epsilon|n+|\epsilon|p)$.
\end{itemize}

We would like to point out that compared to the clustering methods \cite{Wang2021}, our approach has the same computational complexity. 

\section{Numerical Experiments}\label{section5}

This section evaluates  our method against state-of-the-art fusion regularized clustering approaches and other advanced multi-view clustering  techniques. 
Subsection \ref{4.1} lists the experimental setup, Subsection \ref{4.2} provides the experimental results, Subsection \ref{4.3} gives the ablation study, Subsection \ref{4.4} presents the discussions.


\subsection{Experimental Setup}\label{4.1}

\subsubsection{Data Description}

We  describe six commonly used data documented in Table \ref{tab3}:
\begin{itemize}
	\item Authors\footnotemark[1] is  composed of word counts from 841 chapters written by four famous
	English-language authors (Austen, London, Milton, and Shakespeare). 
	\item Lung-discrete\footnotemark[2]  is a biological data which contains 73 observations with 325	features, and each sample belongs to one of 7 classes. 
	\item Dermatology\footnotemark[3]  is a biological data of 6 erythemato-squamous disease diagnoses:  psoriasis, seborrheic dermatitis (SD), lichen planus (LP), pityriasis rosea (PR), chronic dermatitis (CD), and pityriasis rubra pilaris (PRP). It is compiled from clinical assessments and histopathological examinations  of skin samples collected from 358 patients. 
	We perform Z-score standardization on the complete data.
	\item UCI-digits\footnotemark[4] is  an image data that contains 2,000 samples of handwritten digits (``0'' to ``9'').  Each digit is represented by 6 distinct feature sets: profile correlations, Fourier coefficients of the character shapes (FCCS), Karhunen-Lo\`eve coefficients, morphological features, pixel averages in $2\times 3$ windows, and Zernike moments (ZM).  We apply square-root transformations to  FCCS and ZM. 
	\item MSRCV1\footnotemark[5] is  an image data containing 210
	images with 7 classes: tree, car, face, cow, bicycle, building, and airplane. The features of these images are collected from 6 views:
	Census Transform (CENT), CMT, Generalized Search Tree,  HOG, LBP, and Scale Invariant Feature Transform. We take the log of CENT and take the square root of LBP and HOG. 
	Further, we use the Z-score method to standardize the entire data.
	\footnotetext[1]{https://github.com/DataSlingers/clustRviz/tree/master/data}
	\footnotetext[2]{https://jundongl.github.io/scikit-feature/datasets.html}
	\footnotetext[3]{http://archive.ics.uci.edu/ml/datasets/Dermatology}
	\footnotetext[4]{http://archive.ics.uci.edu/ml/datasets/Multiple+Features}
	\footnotetext[5]{https://github.com/Jeaninezpp/multi-view-datasets}
	\item Multi-omics\footnotemark[6] is a  biological data on 343 breast cancer tumor samples from TCGA, which can be divided into Luminal, Basal-like, and HER2-enriched subtype based on PAM50 mRNA classification\footnotemark[7]. It includes 4  views:  RNA gene expression profiles, DNA methylation data, miRNA expression levels, and reverse phase protein array data. All data preprocessing follows \cite{Wang2021}.
	\footnotetext[6]{https://github.com/ttriche/bayesCC}
	\footnotetext[7]{https://pmc.ncbi.nlm.nih.gov/articles/PMC3465532/}
\end{itemize}
\begin{table}[width=.9\linewidth,cols=5,pos=t]
	\caption{Details of real data.}
	\label{tab3}
	\begin{tabular*}{\tblwidth}{@{} LLLLL@{} }
		\toprule
		Data & Clusters &Features & Observations & Features of per view (total views)\\ 
		\midrule
		Authors & 4 & 69 & 841 & 69 (1) \\ 
		Lung-discrete& 7  & 325 & 73 & 325 (1) \\ 
		Dermatology & 6 & 34 & 358 & 12, 22 (2) \\
		UCI-digits & 10 & 649 & 2000 & 216, 76, 64, 6, 240, 47 (6) \\
		MSRCV1  & 7 & 2428 &210  & 1302, 48, 512, 100, 256, 210 (6) \\
		Multi-omics & 3  & 1813 & 343 & 645, 423, 574, 171 (4) \\ 
		\bottomrule
	\end{tabular*}
\end{table}

\subsubsection{Comparing Methods}

To highlight the merits of our approach,  we conduct a series of methodological comparisons. 
\begin{itemize} 
		\item $K$-means:  $K$-means clustering. It is a baseline method, whose experimental result is obtained by executing the built-in function kmeans in MATLAB. 
		\item SCC \cite{Wang2018}:   Sparse convex fusion regularized clustering  with inter-group Lasso penalty in \eqref{YY2}.
		\item SGLCC \cite{Chen2020}: Convex fusion regularized clustering  with sparse group Lasso penalty.
		\item iGecco+ \cite{Wang2021}: Multi-view fusion regularized clustering  with inter-group Lasso penalty in \eqref{YY3}. 
		\item SCMvFS   \cite{Cao2024}: Structure learning with consensus label information and heterogeneous graphs for multi-view feature selection.
	The model employs smooth fusion regularization terms in the learning of the clustering index matrix. 
	\item CDMvFS \cite{CAO202412}:  Multi-view unsupervised feature selection that enforces mutual orthogonality among diverse graphs while  identifying a consensus partition.
	\item WMVEC  \cite{LIU202411}: Adaptive weighted multi-view evidential clustering employing belief function theory and feature preference to address uncertainty/imprecision in cluster assignments and feature irrelevance/redundancy.
\end{itemize}
\begin{figure*}[t]
	\centering
	\hspace*{-0.3cm} 
	\subfloat[]{
		\centering
		\includegraphics[width=0.32\textwidth]{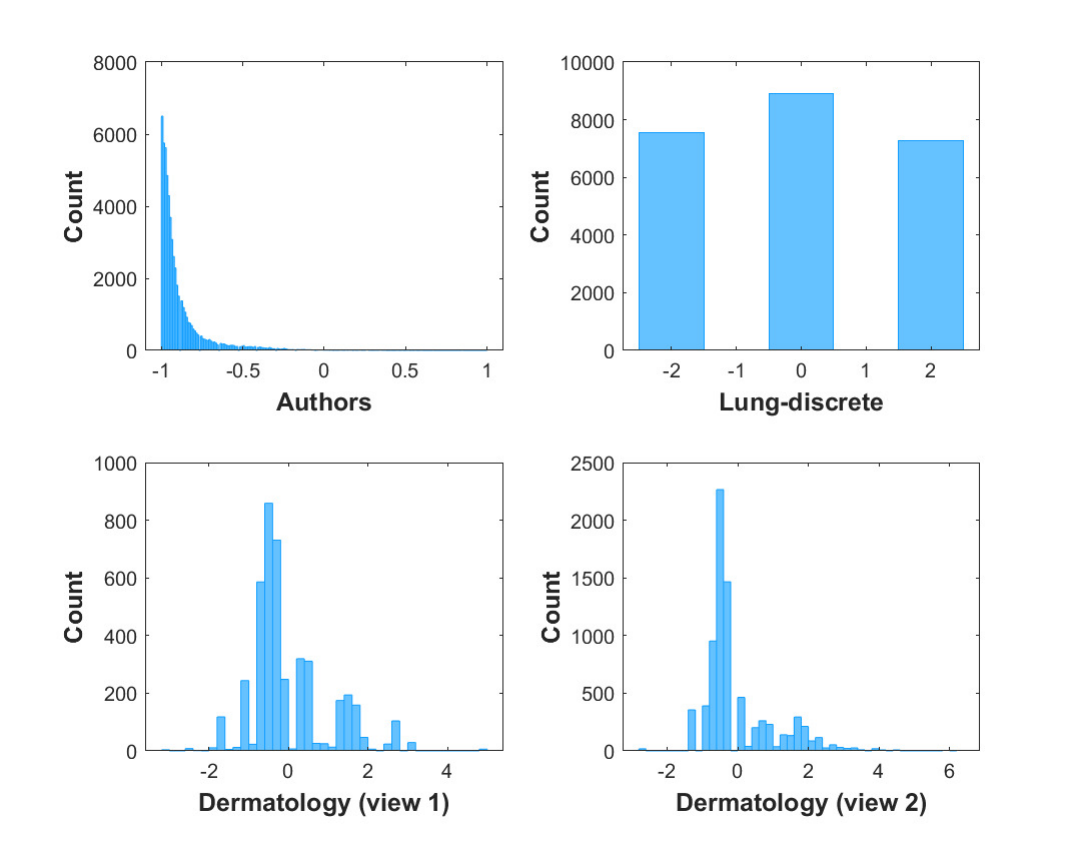}
		\label{fig:subfig7}
	}
	\subfloat[]{
		\centering
		\includegraphics[width=0.5\textwidth]{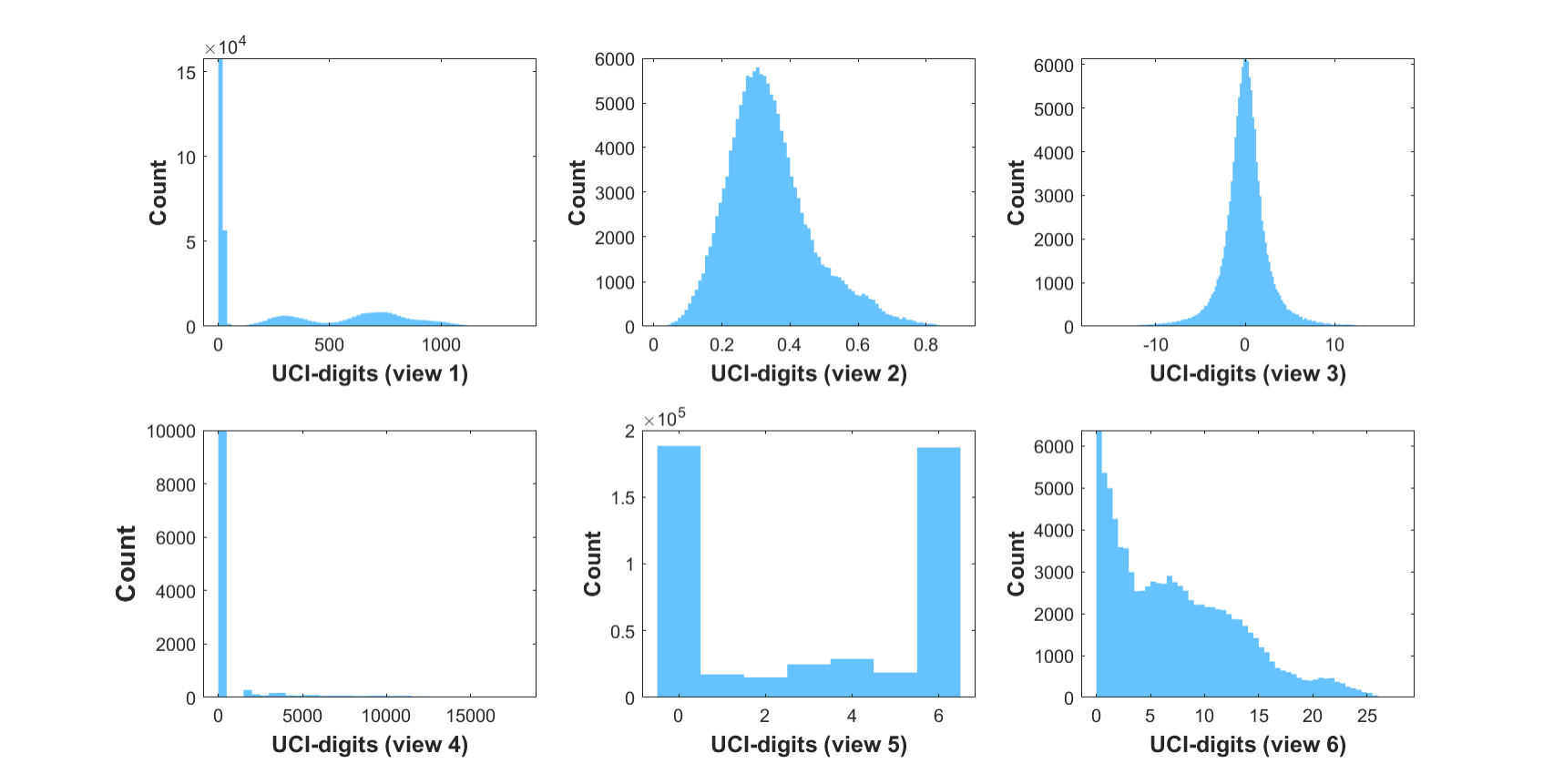}
		\label{fig:subfig8}
	}\\	\vspace{-0.7cm}
	\subfloat[]{
		\centering
		\includegraphics[width=0.32\textwidth]{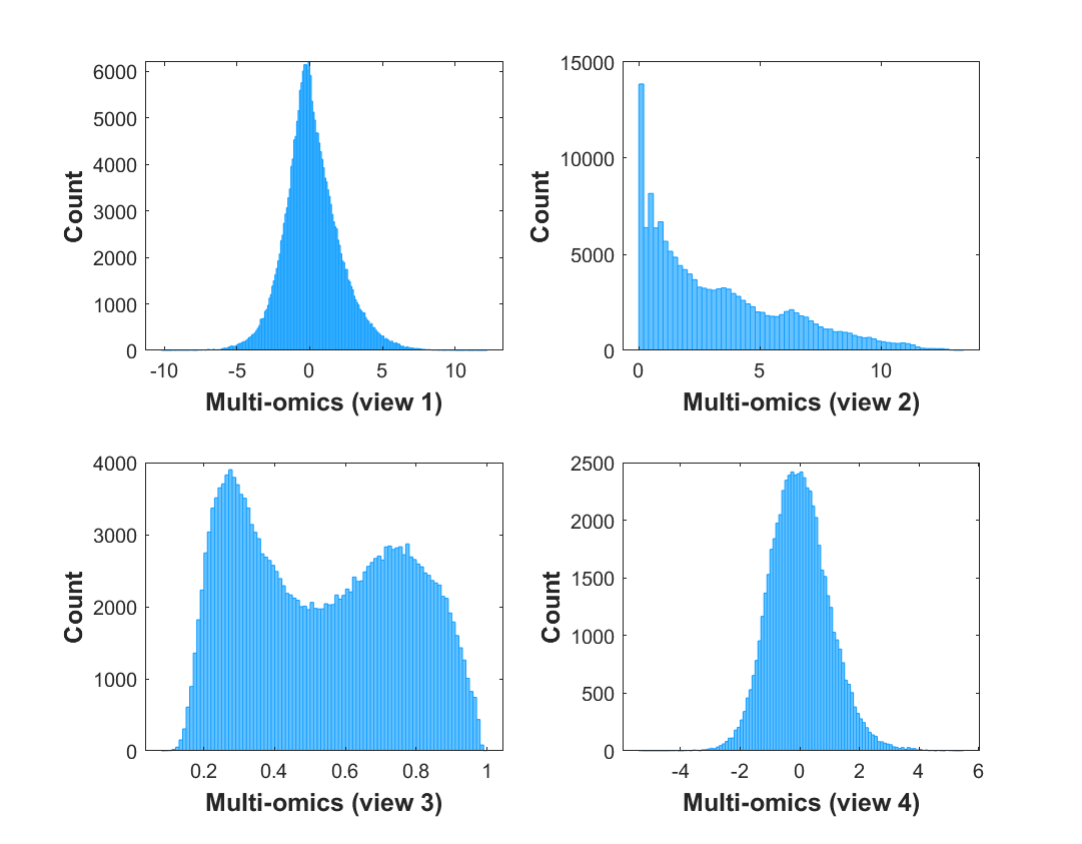}
		\label{fig:subfig9}
	}
	\subfloat[]{
		\centering
		\includegraphics[width=0.5\textwidth]{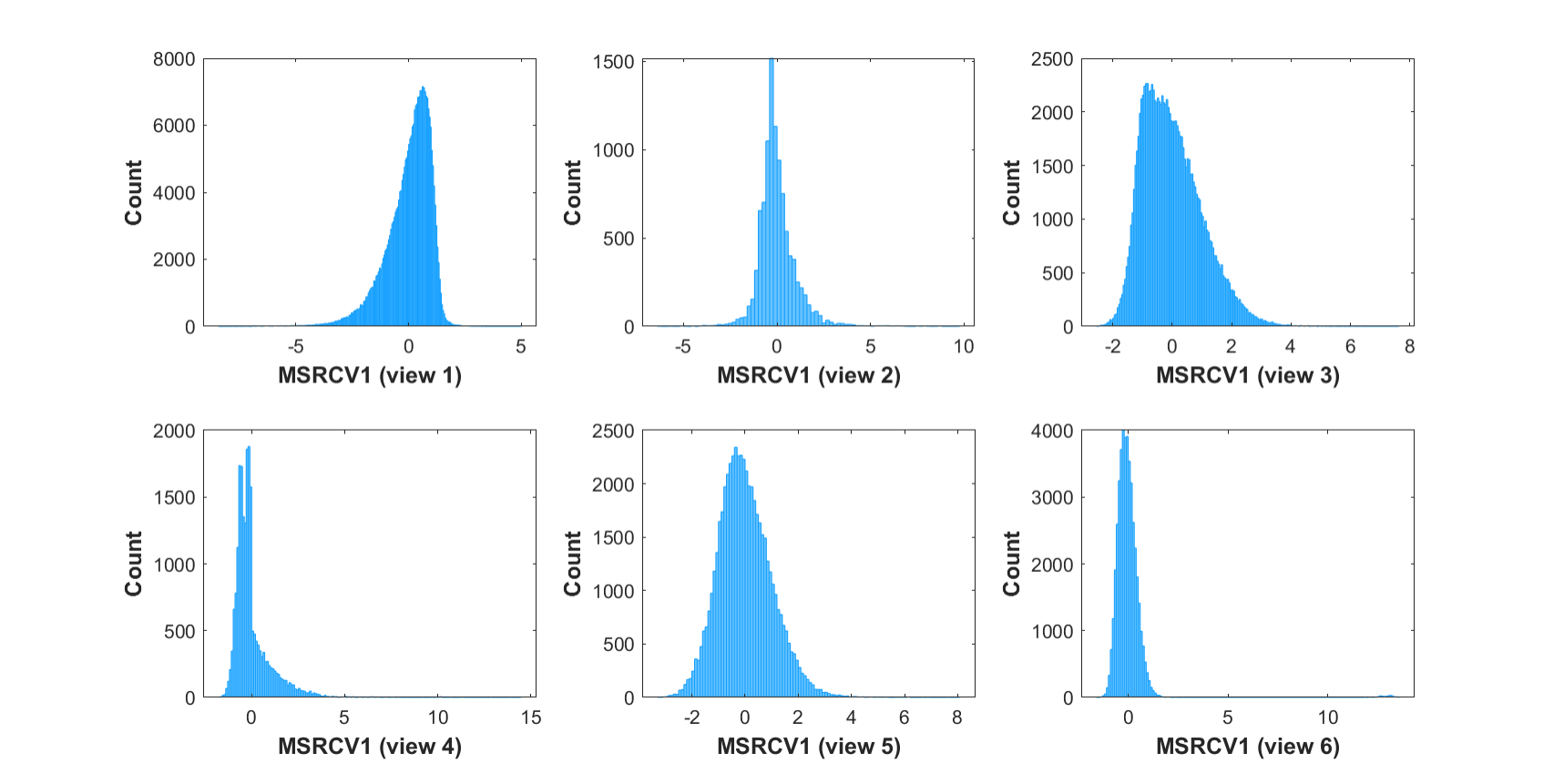}
		\label{fig:subfig10}
	}
	\caption{Histograms of real data showing heterogeneous characteristics.  Authors is highly-skewed. Lung-discrete and Dermatology appear discrete Gaussian.  UCI-digits shows view-specific heterogeneity (skewed views 1/4/6, Gaussian views 2/3, and count data view 5). MSRCV1 exhibits diverse and potentially non-Gaussian across views.  Multi-omics displays the mixture distribution characteristics, including Gaussian gene expression (view 1) and protein (view 4), right-skewed miRNA (view 2), and proportion-valued methylation (view 3).}
	\label{fig4}
\end{figure*}

Unless specified otherwise, the reported $K$-means, SCMvFS, and CDMvFS results are averaged across 30 executions. The other methods are not stochastic.
\begin{table*}[width=.9\linewidth,cols=10,pos=t]
	\caption{Comparisons of ACC, NMI, ARI and FMI for different clustering methods on  real data.}
	\label{tab4}
	\begin{tabular*}{\tblwidth}{@{} LLLLLLLLLL@{} }
		\toprule
		Data&Methods &$K$-means & SCC & SGLCC & iGecco+  &SCMvFS& CDMvFS&WMVEC  &Our \\
		\midrule
		\multirow{4}{*}{\makecell{Authors}}	
		&ACC $\uparrow$ &0.7915 & 0.9203  & \underline{0.9952}  &0.9906  & 0.8950&0.8067 & 0.9620   &\bfseries 0.9988\\
		&NMI $\uparrow$&0.7118 &0.8162 & \underline{0.9768} &  0.9544&0.6967 & 0.7055 &0.8669  	&\bfseries 0.9936\\
		&ARI \;$\uparrow$&0.7221 & 0.8941  & \underline{0.9858}  &0.9698& 0.7258 &0.7212 & 0.8943    &\bfseries 0.9969\\
		&FMI \,$\uparrow$&0.8051 &0.9291 & \underline{0.9902} &  0.9792  &0.8096 & 0.8050 &0.9269 	&\bfseries 0.9979\\
		\midrule
		\multirow{4}{*}{\makecell{Lung-discrete}} 
		&ACC $\uparrow$ &0.6521&0.6575 &\bfseries 0.8630 &0.7134 & 0.7324 &0.7256&\underline{0.7397}   &\bfseries 0.8630  \\
		&NMI $\uparrow$ &0.6258&  0.5563& \underline{0.7777}& 0.6399 &0.7059 &0.6723&0.7399	& \bfseries 0.7812  \\&		
		ARI \;$\uparrow$ &0.4778 &0.4831 &\underline{0.7236} &0.5516 & 0.6288 &0.5723& 0.6535  &\bfseries 0.7345  \\
		&FMI \,$\uparrow$&0.5703 &  0.6266&  \underline{0.7724}& 0.6668&0.6923&0.6459&0.7161	& \bfseries 0.7809  \\
		\midrule
		\multirow{4}{*}{\makecell{Dermatology}} 
		&ACC $\uparrow$ &0.6698&0.8601  &0.8547 &\underline{0.8603}  &0.7989 &0.7821 &0.8268 &	 \bfseries 0.8631\\	
		&NMI $\uparrow$&0.7820& 0.8619& 0.8540 &  \underline{0.8674}&0.8392&0.8347&0.8445 &\bfseries 0.8769\\
		&ARI \;$\uparrow$ &0.6202 &0.8469  &0.8377 &\underline{0.8615} &0.7604 &0.7456 &0.8293  &	 \bfseries 0.8694\\	
		&FMI \,$\uparrow$&0.7083& 0.8806& 0.8777 &  \underline{0.8967} &0.8125&0.8021&0.8713&\bfseries 0.9028\\
		\midrule
		\multirow{4}{*}{\makecell{UCI-digits}} 
		&ACC $\uparrow$ &0.5025 &0.5150  &0.5145 &0.6710&\underline{0.7686} &\bfseries0.8151  &0.6920  &\underline{0.7686}\\	
		&NMI $\uparrow$&0.5630 & 0.6160& 0.6149 &  0.7222 &0.7371 &\underline{0.7570}&0.6267 &\bfseries 0.7634\\
		&ARI \;$\uparrow$ &0.4213&0.4827 &0.4829 &0.6206 &0.6639 &\underline{0.6871} &0.5232  &	  \bfseries0.6965 \\	
		&FMI \,$\uparrow$&0.4862& 0.5981 & 0.5992 &  0.6952&0.7000	&\underline{0.7197}&0.5711 &\bfseries 0.7355\\
		\midrule
		\multirow{4}{*}{\makecell{MSRCV1 }} 
		&ACC $\uparrow$ &0.6251 &0.6238  &0.4667  &0.7270 &0.7267&\bfseries0.7975&0.6524   & \underline{0.7762}\\	
		&NMI $\uparrow$ &0.5824& 0.6740& 0.4606 &  0.7067&0.6676&\underline{0.7266}&0.6541 &\bfseries 0.7504\\
		&ARI \;$\uparrow$ &0.4619 &0.5299  &0.2955  &0.5853  &0.5825&\underline{0.6559}&0.5541  & \bfseries 0.6684\\	
		&FMI \,$\uparrow$&0.5531& 0.6230& 0.4681 &  0.6723&0.6527&\underline{0.7078}&0.6893 &\bfseries 0.7305\\
		\midrule
		\multirow{4}{*}{\makecell{Multi-omics}} &	
		ACC $\uparrow$ &0.6650 &0.8630  &0.8659  &\underline{0.9329}&0.9246&0.9162 &0.8338   &	 \bfseries 0.9504\\	
		&NMI $\uparrow$&0.4932& 0.5296& 0.5292 &  0.7019&\underline{0.7086}	&0.6391&0.4676 &\bfseries 0.7575\\
		&ARI \;$\uparrow$ &0.4304 &0.6185  &0.6194  &\underline{0.8066} &0.7875&0.7723 &0.6726   &	 \bfseries 0.8567\\	
		&FMI \,$\uparrow$&0.6875& 0.8549& 0.8529 &  \underline{0.9149}&0.8966	&0.8855&0.8460 &\bfseries 0.9350\\
		\bottomrule
	\end{tabular*}
\end{table*}
\begin{figure}[h]
	\centering
	\captionsetup[subfigure]{font=scriptsize, textfont=scriptsize, labelfont=scriptsize}
	\hspace*{-0.3cm} 
	\subfloat[(a1) K-means]{
		\centering
		\includegraphics[width=0.122\textwidth]{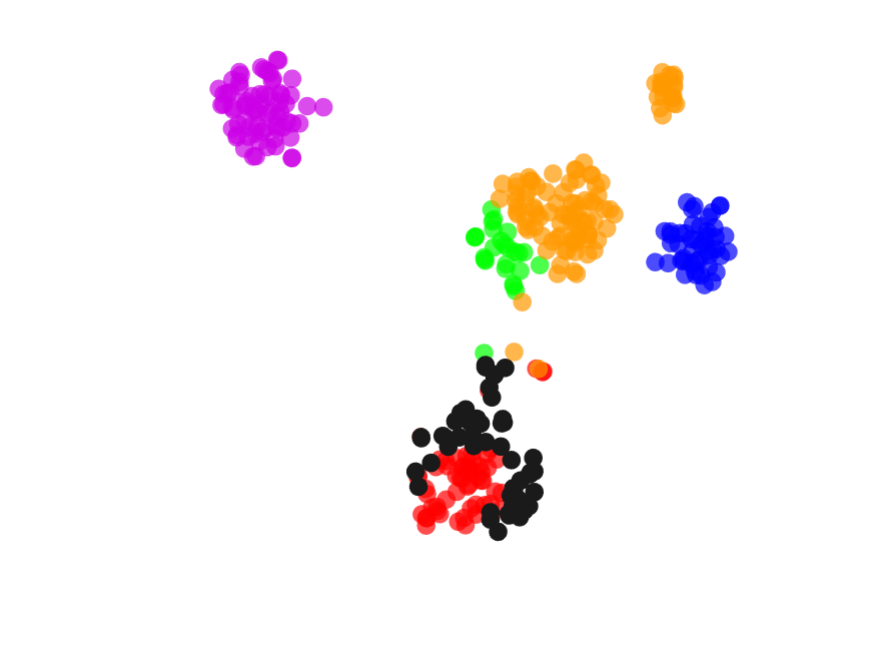}
		\label{fig:subfig21}
	}	\hspace*{-0.3cm} 
	\subfloat[(a2) SCC]{
		\centering
		\includegraphics[width=0.122\textwidth]{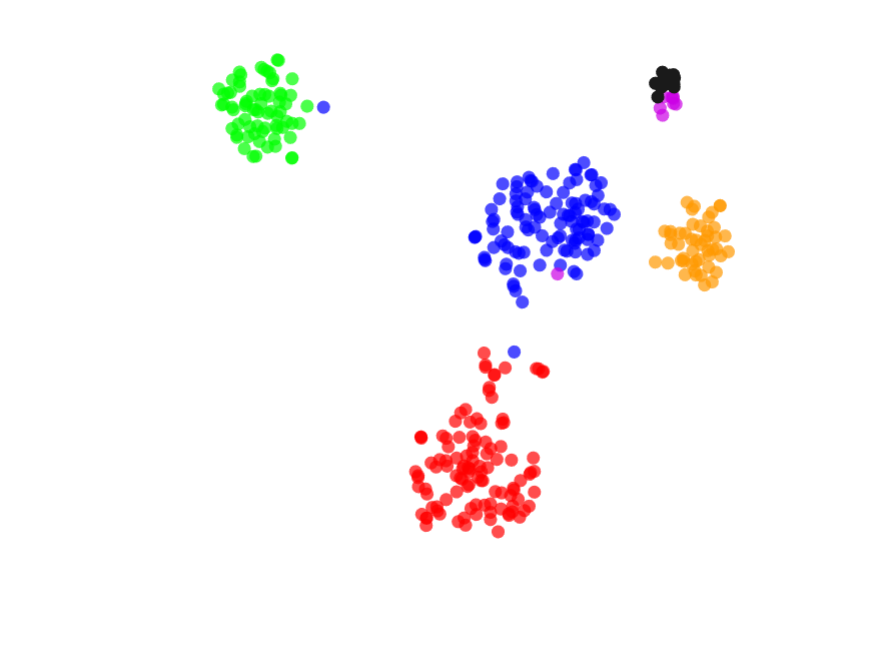}
		\label{fig:subfig25}
	}	\hspace*{-0.3cm} 
	\subfloat[(a3) SGLCC]{
		\centering
		\includegraphics[width=0.122\textwidth]{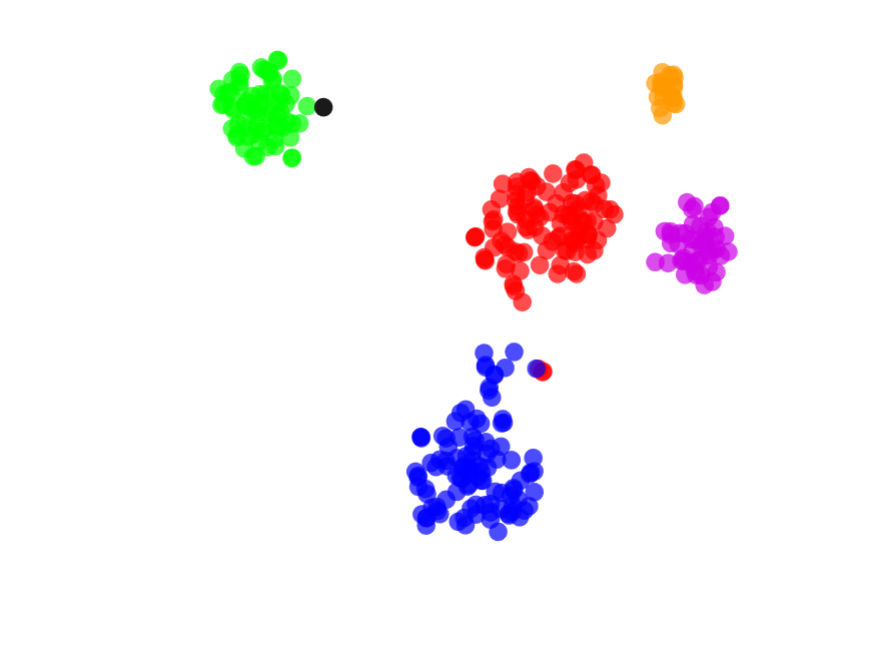}
		\label{fig:subfig26}
	}	\hspace*{-0.3cm} 
	\subfloat[(a4) iGecco+]{
		\centering
		\includegraphics[width=0.122\textwidth]{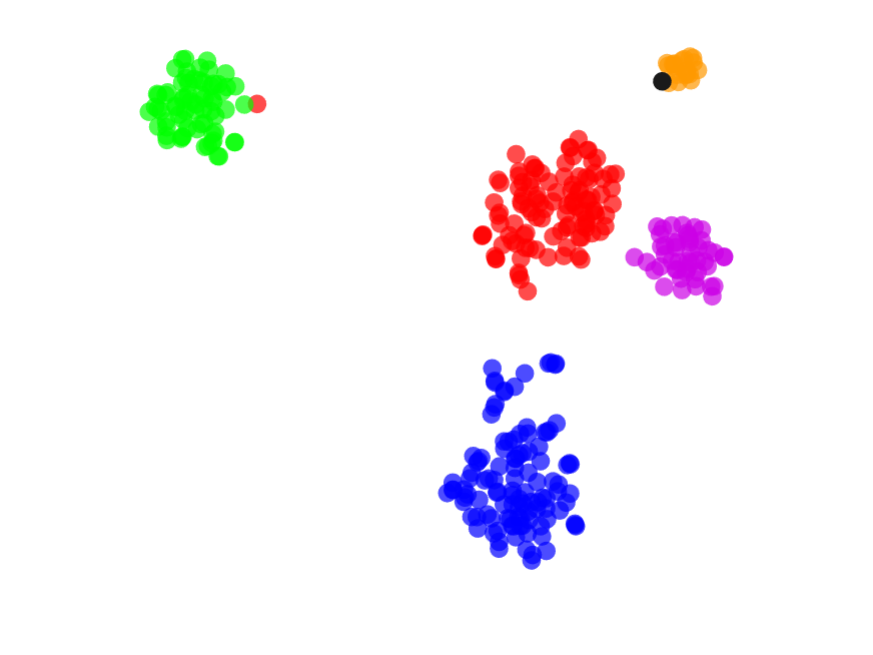}
		\label{fig:subfig27}
	}	\hspace*{-0.3cm} 
	\subfloat[(a5) SCMvFS]{
		\centering
		\includegraphics[width=0.122\textwidth]{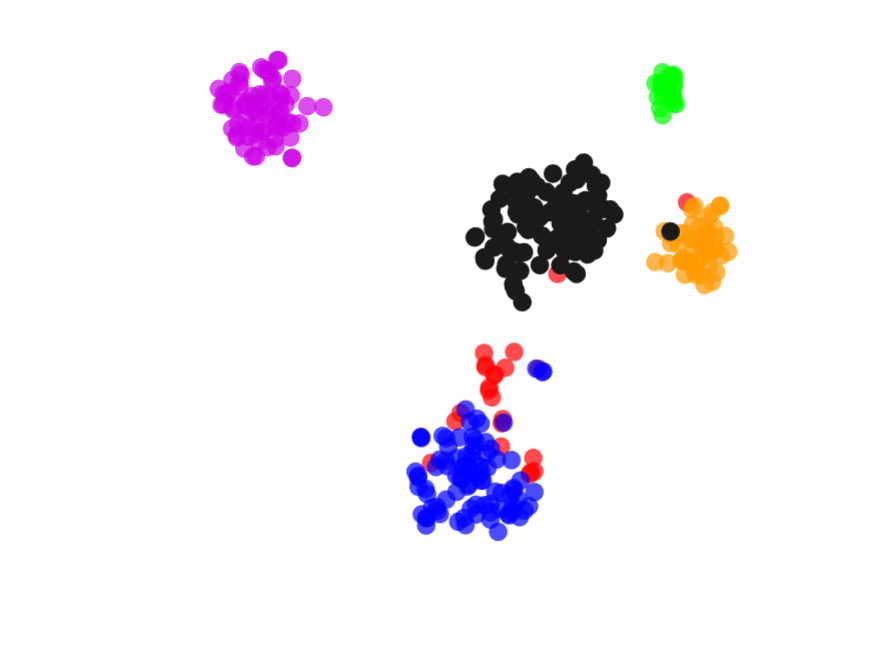}
		\label{fig:subfig23}
	}	\hspace*{-0.3cm} 
	\subfloat[(a6) CDMvFS]{
		\centering
		\includegraphics[width=0.122\textwidth]{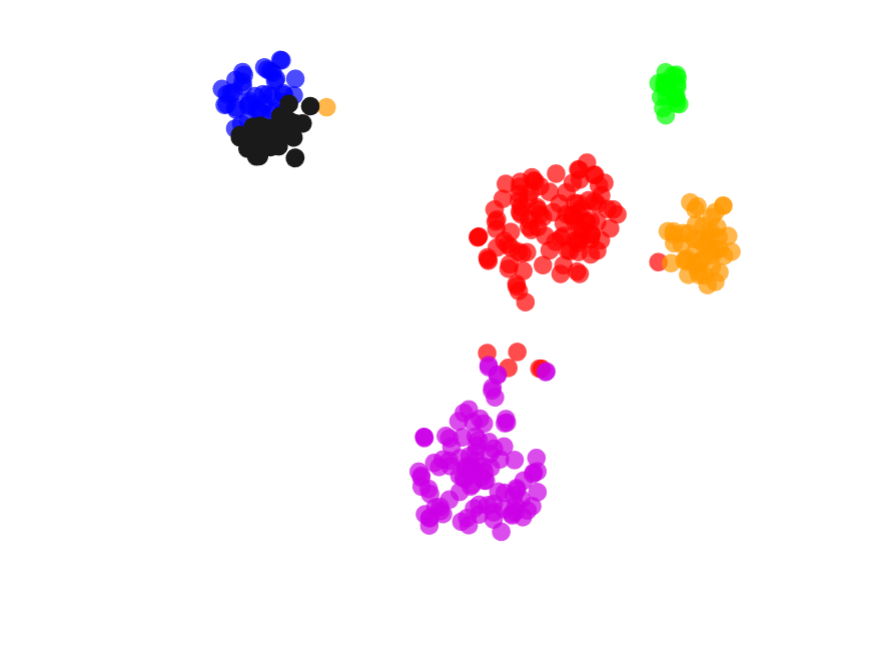}
		\label{fig:subfig22}
	}	\hspace*{-0.3cm} 
	\subfloat[(a7) WMVEC]{
		\centering
		\includegraphics[width=0.122\textwidth]{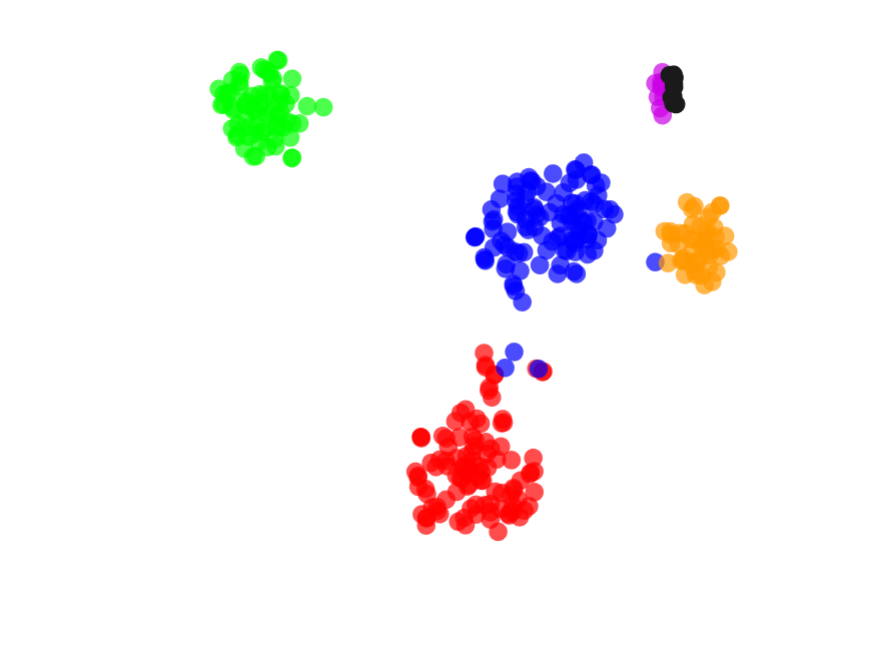}
		\label{fig:subfig24}
	}	\hspace*{-0.3cm} 
	\subfloat[(a8) Our]{
		\centering
		\includegraphics[width=0.122\textwidth]{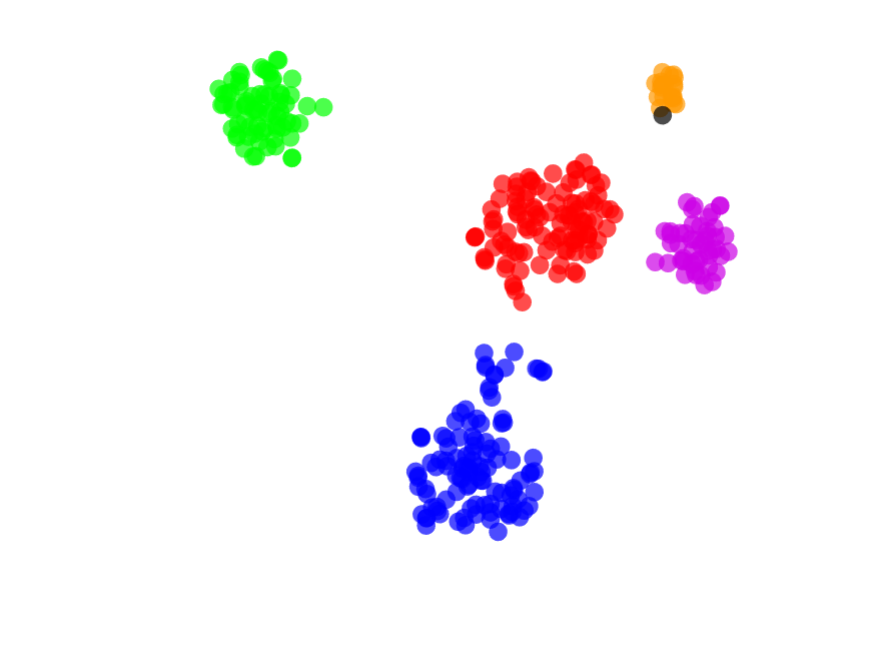}
		\label{fig:subfig28}
	}
	\\
	\hspace*{-0.3cm} 
	\subfloat[(b1) K-means]{
		\centering
		\includegraphics[width=0.122\textwidth]{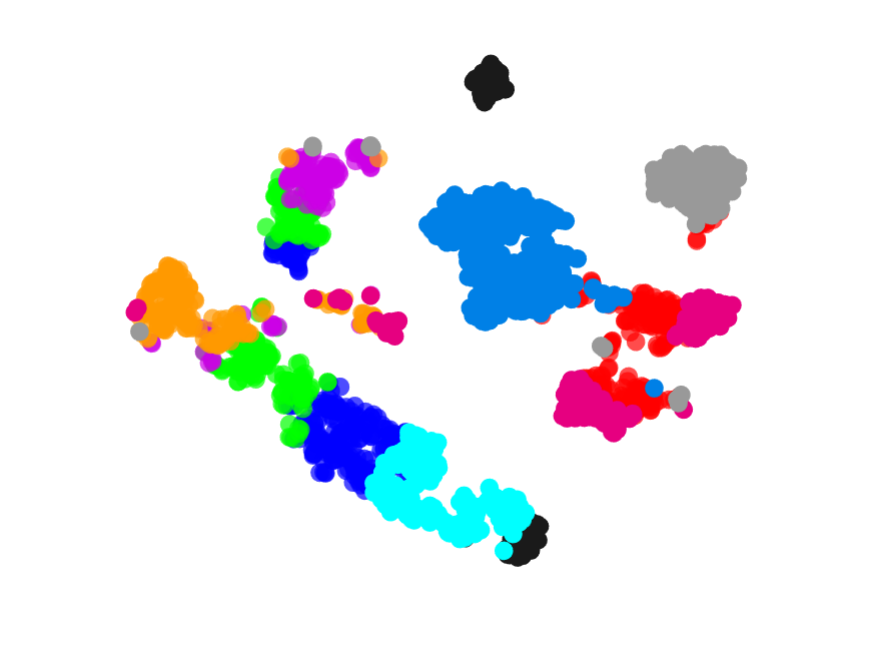}
		\label{fig:subfig61}
	}	\hspace*{-0.3cm} 
	\subfloat[(b2) SCC]{
		\centering
		\includegraphics[width=0.122\textwidth]{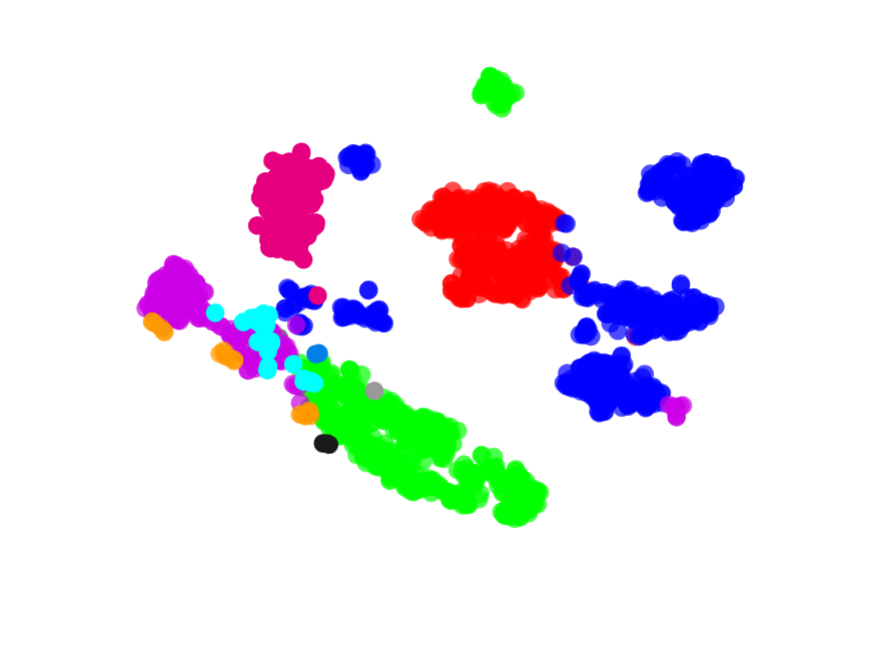}
		\label{fig:subfig62}
	}	\hspace*{-0.3cm} 
	\subfloat[(b3) SGLCC]{
		\centering
		\includegraphics[width=0.122\textwidth]{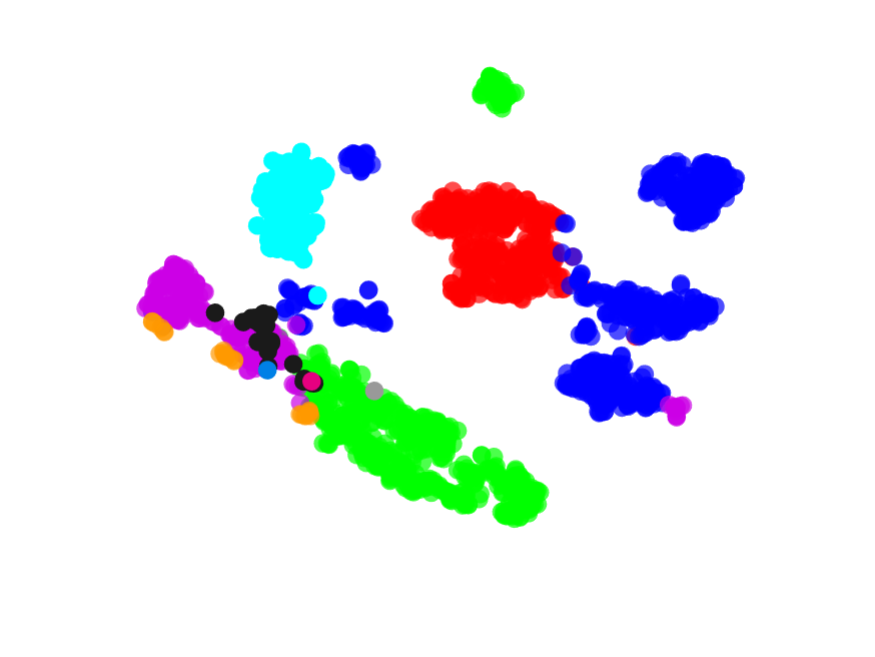}
		\label{fig:subfig63}
	}	\hspace*{-0.3cm} 
	\subfloat[(b4) SCMvFS]{
		\centering
		\includegraphics[width=0.122\textwidth]{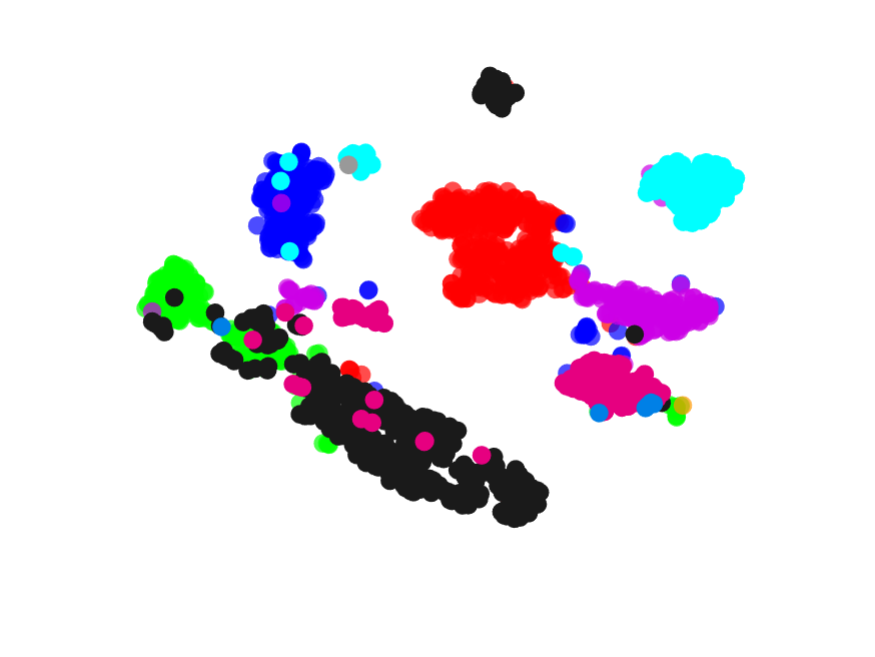}
		\label{fig:subfig64}
	}	\hspace*{-0.3cm} 
	\subfloat[(b5) SCMvFS]{
		\centering
		\includegraphics[width=0.122\textwidth]{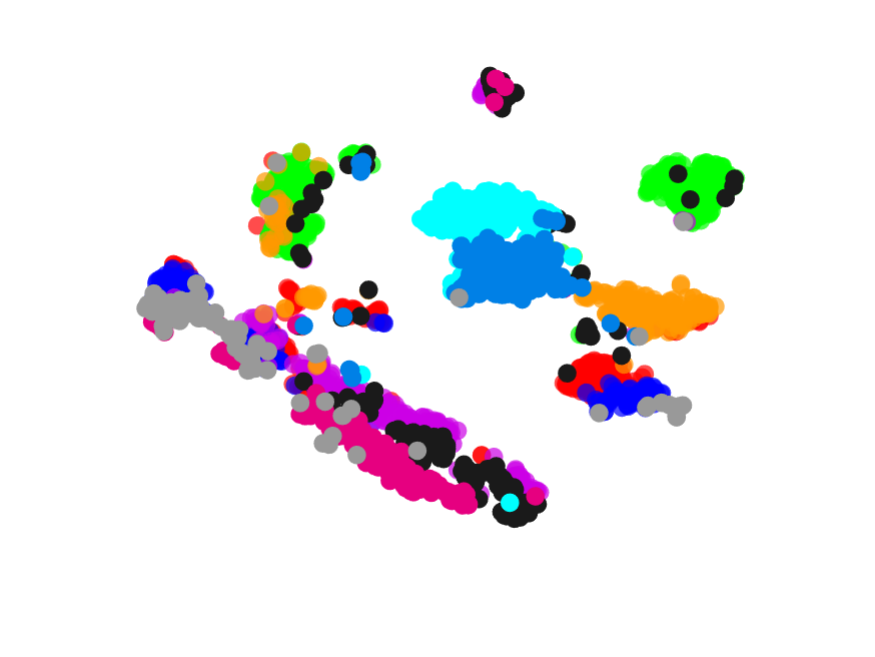}
		\label{fig:subfig65}
	}	\hspace*{-0.3cm} 
	\subfloat[(b6) CDMvFS]{
		\centering
		\includegraphics[width=0.122\textwidth]{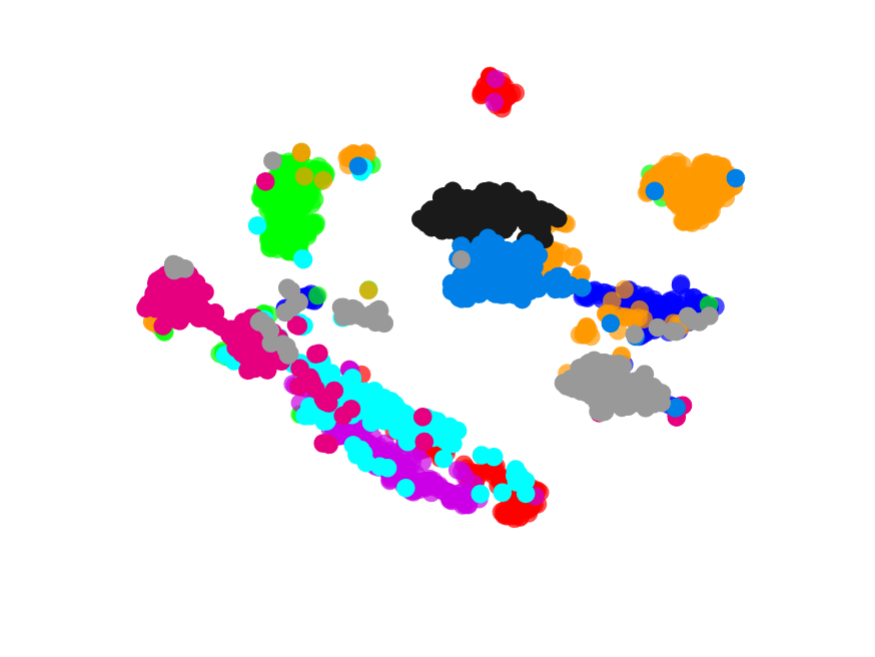}
		\label{fig:subfig66}
	}	\hspace*{-0.3cm} 
	\subfloat[(b7) WMVEC]{
		\centering
		\includegraphics[width=0.122\textwidth]{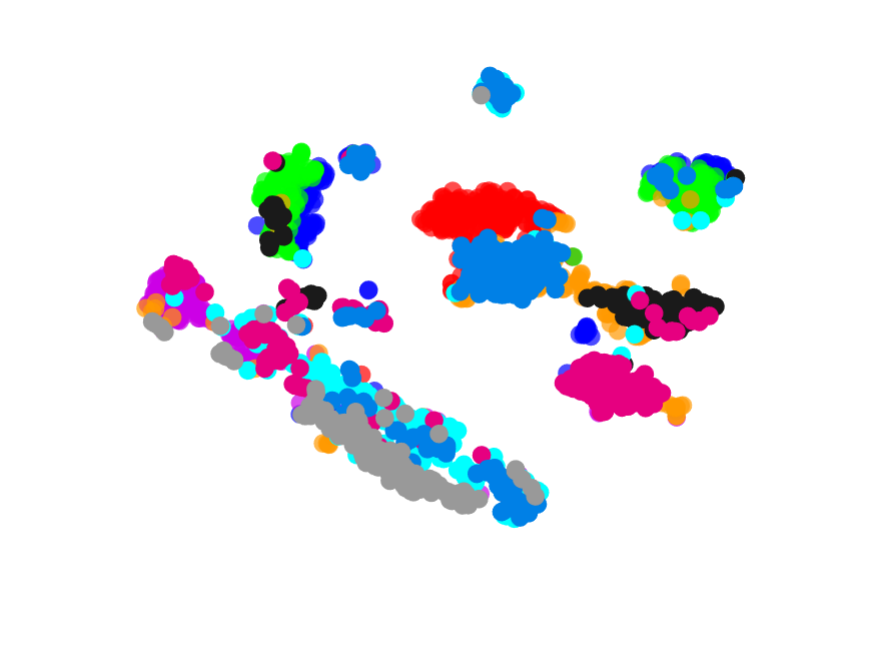}
		\label{fig:subfig67}
	}	\hspace*{-0.3cm} 
	\subfloat[(b8) Our]{
		\centering
		\includegraphics[width=0.122\textwidth]{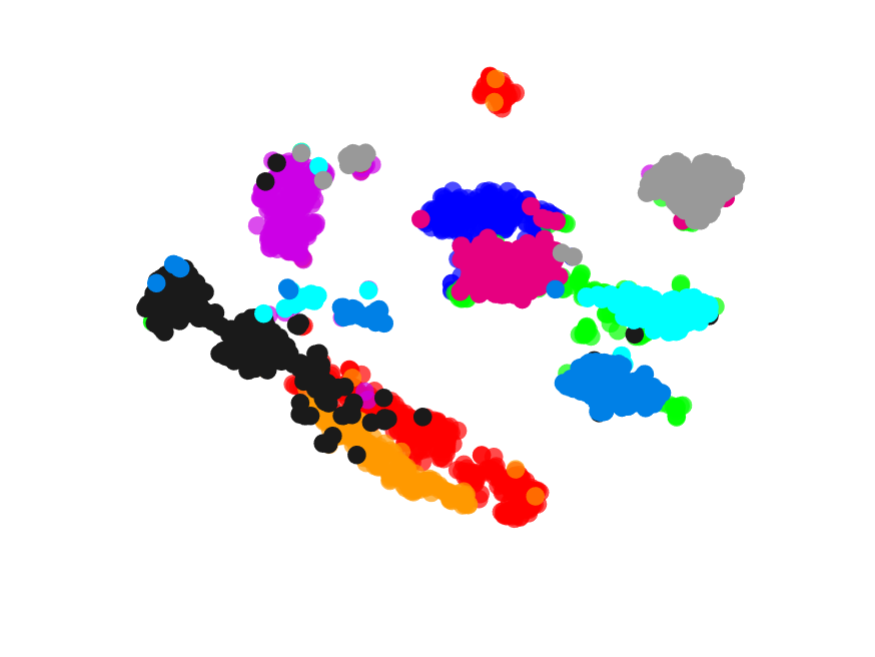}
		\label{fig:subfig68}
	}
	\\
	\hspace*{-0.3cm} 
	\subfloat[(c1) K-means]{
		\centering
		\includegraphics[width=0.122\textwidth]{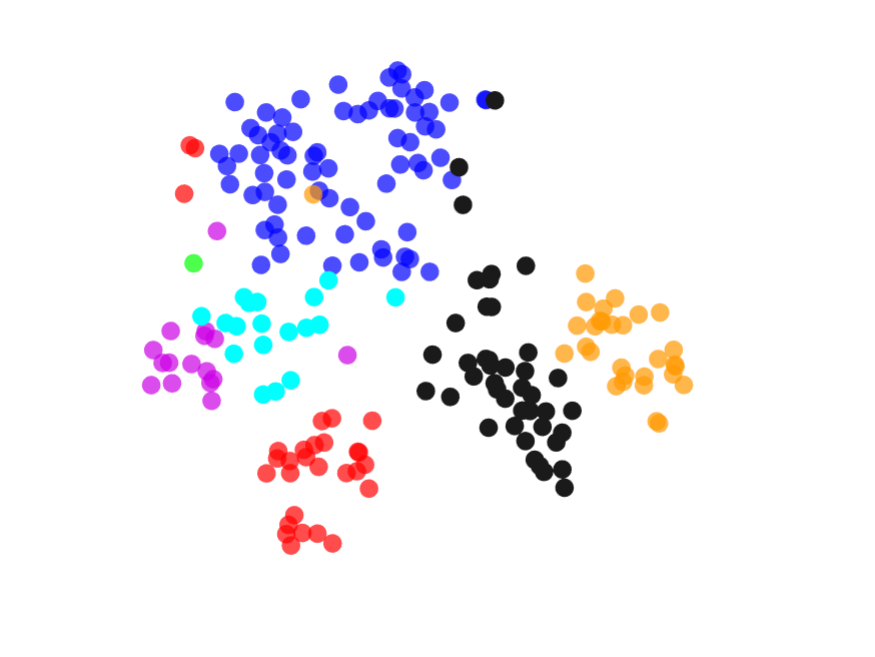}
		\label{fig:subfig51}
	}	\hspace*{-0.3cm} 
	\subfloat[(c2) SCC]{
		\centering
		\includegraphics[width=0.122\textwidth]{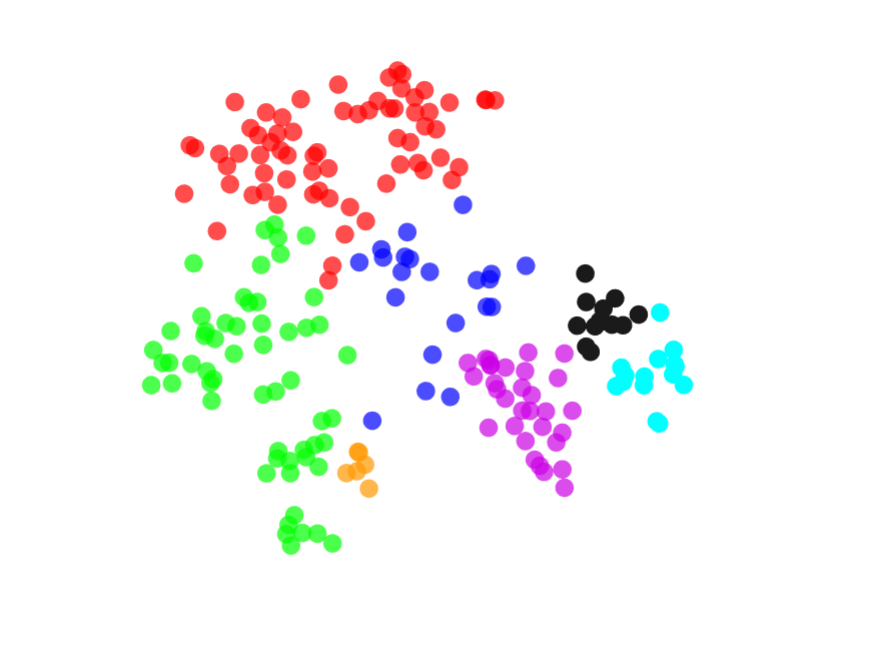}
		\label{fig:subfig52}
	}	\hspace*{-0.3cm} 
	\subfloat[(c3) SGLCC]{
		\centering
		\includegraphics[width=0.122\textwidth]{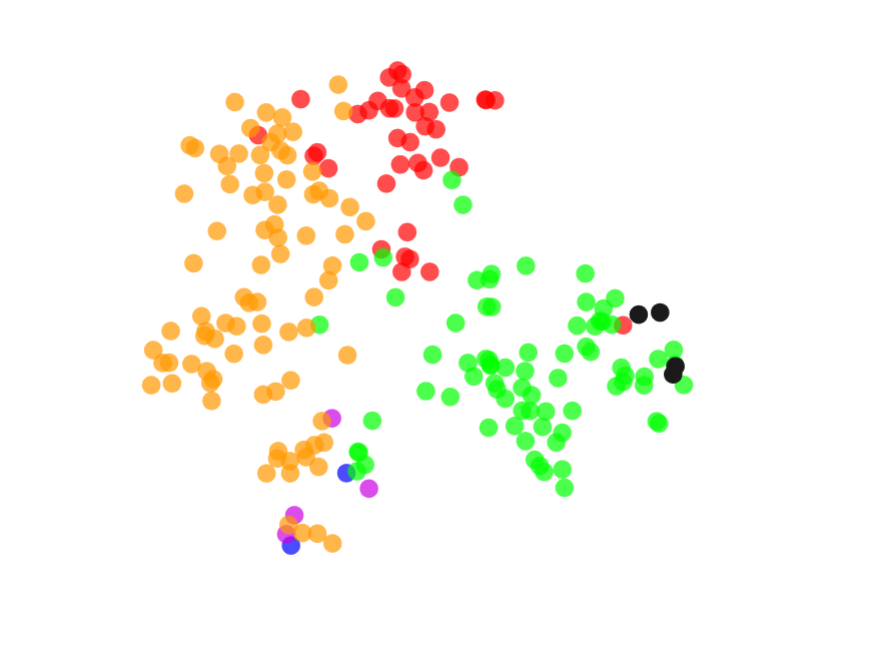}
		\label{fig:subfig53}
	}	\hspace*{-0.3cm} 
	\subfloat[(c4) SCMvFS]{
		\centering
		\includegraphics[width=0.122\textwidth]{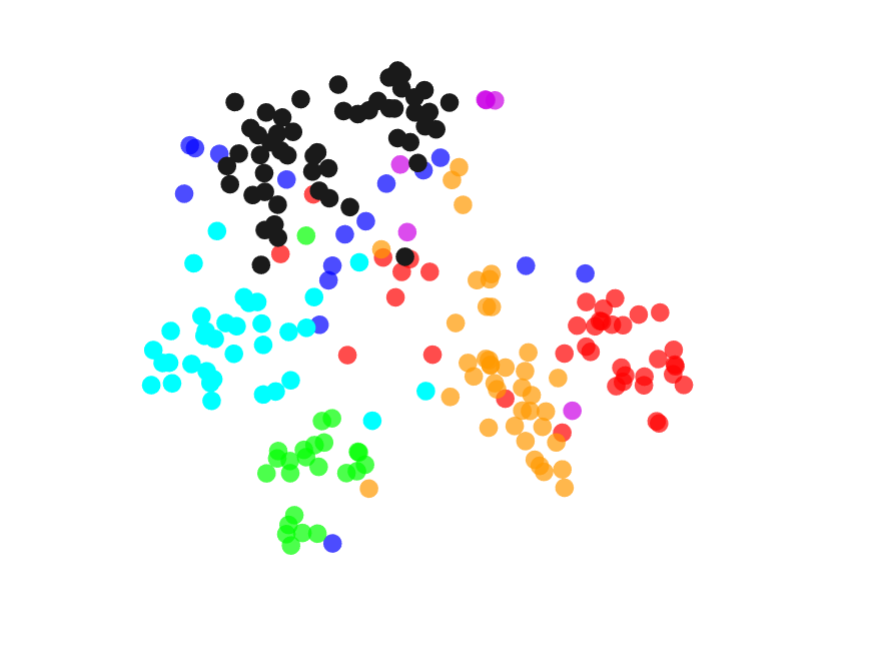}
		\label{fig:subfig54}
	}	\hspace*{-0.3cm} 
	\subfloat[(c5) SCMvFS]{
		\centering
		\includegraphics[width=0.122\textwidth]{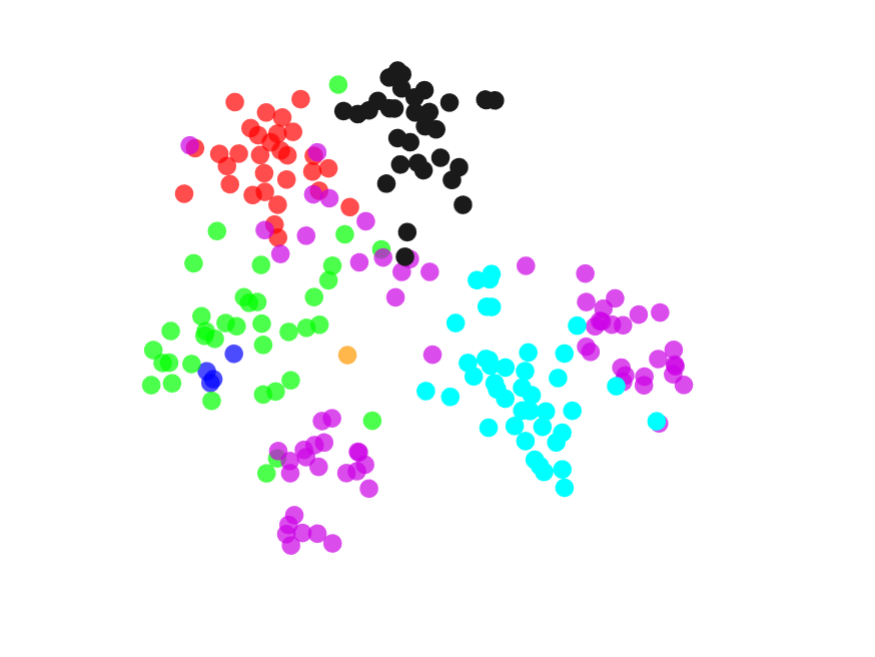}
		\label{fig:subfig55}
	}	\hspace*{-0.3cm} 
	\subfloat[(c6) CDMvFS]{
		\centering
		\includegraphics[width=0.122\textwidth]{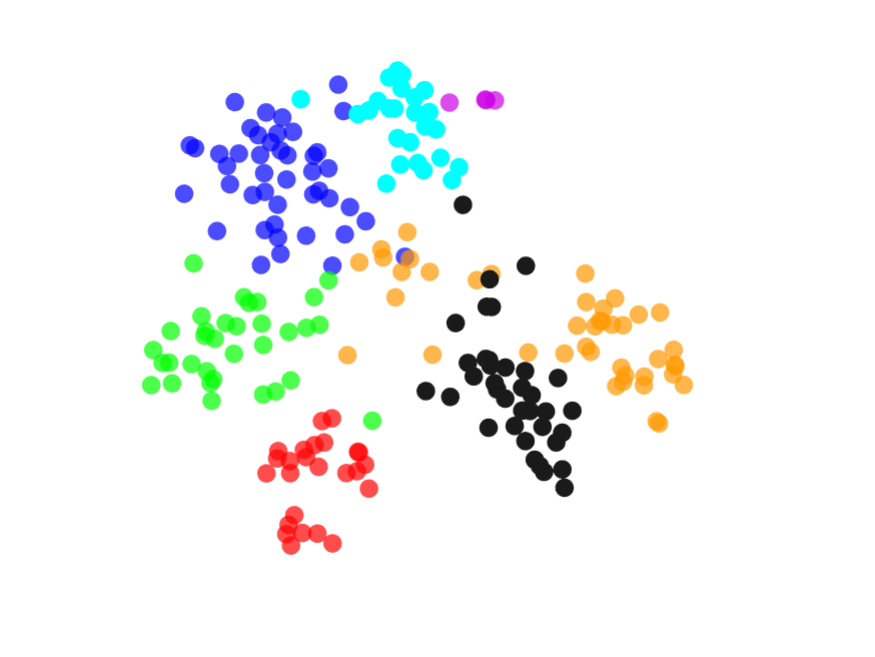}
		\label{fig:subfig56}
	}	\hspace*{-0.3cm} 
	\subfloat[(c7) WMVEC]{
		\centering
		\includegraphics[width=0.122\textwidth]{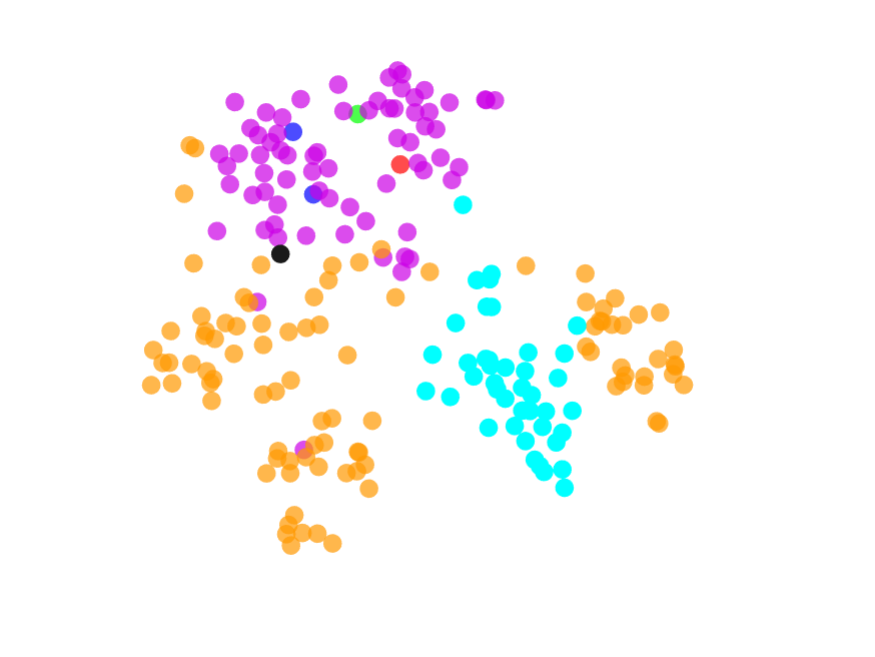}
		\label{fig:subfig57}
	}	\hspace*{-0.3cm} 
	\subfloat[(c8) Our]{
		\centering
		\includegraphics[width=0.122\textwidth]{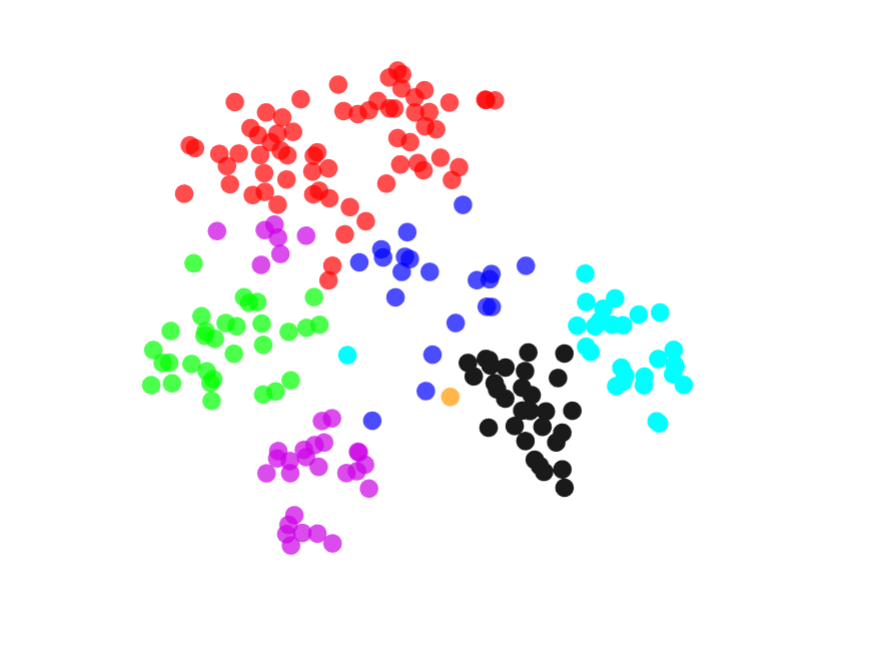}
		\label{fig:subfig58}
	}
	\caption{Visual comparison of clustering results from 8 methods on  Dermatology (a1-a8), UCI-digits (b1-b8), and MSRCV1 (c1-c8).}
	\label{fig2}
\end{figure}
\begin{table*}[width=.9\linewidth,cols=9,pos=t]
	\caption{Ablation study of clustering performance and feature selection capability.}
	\label{tab6}
	\begin{tabular*}{\tblwidth}{@{} LLLLLLLLL@{} }
		\toprule
		Data&$\|\cdot\|_{2,0}$& $\|\cdot\|_{0}$ &ACC $\uparrow$& NMI $\uparrow$&ARI $\uparrow$& FMI $\uparrow$&Selected feature number & $|T_{cur}\cap T_{our}|$ \\
		\midrule
		\multirow{4}{*}{\makecell{Authors}}	&\checkmark &\checkmark&\bfseries0.9988 & \bfseries0.9936&\bfseries0.9969 & \bfseries0.9979& 3  & 3  \\
		&\ding{53}&\checkmark&0.7907&0.6817&0.6741&0.8122&9&3 \\
		&\checkmark &\ding{53}&\bfseries0.9988&\underline{0.9935}&0.\underline{9960}&\underline{0.9972}&1&1 \\
		&\ding{53}&\ding{53}&\underline{0.9905}&0.9540&0.9714&0.9803&69&3 \\
		\midrule
		\multirow{4}{*}{\makecell{Lung-discrete}} &\checkmark &\checkmark&\bfseries0.8630&\bfseries0.7812&\bfseries0.7345&\bfseries0.7809&73&73  \\
		&\ding{53}&\checkmark&\underline{0.8356}& \underline{0.7413}&\underline{0.6912}& \underline{0.7460} &161&47 \\
		&\checkmark &\ding{53}&0.8219&0.7311&0.6710&0.7325&8&5   \\
		&\ding{53}&\ding{53}&0.6986&0.6415&0.5587&0.6639&325&73  \\
		\midrule
		\multirow{4}{*}{\makecell{Dermatology}}	&\checkmark &\checkmark&\bfseries0.8631 & \bfseries0.8769&\bfseries0.8694 & \bfseries0.9028& 3  & 3  \\
		&\ding{53}&\checkmark&\underline{0.8603}&\underline{0.8674}&\underline{0.8615}&\underline{0.8967} &3 &1 \\
		&\checkmark &\ding{53}&\underline{0.8603}&\underline{0.8674}&\underline{0.8615}&\underline{0.8967} &1 &1  \\
		&\ding{53}&\ding{53}&\underline{0.8603}&\underline{0.8674}&\underline{0.8615}&\underline{0.8967} &34 &3 \\
		\midrule
		\multirow{4}{*}{\makecell{UCI-digits}} &\checkmark &\checkmark&\bfseries0.7686&\bfseries0.7634&\bfseries0.6965&\bfseries0.7355&9 &9  \\
		&\ding{53}&\checkmark&\underline{0.7120}& \underline{0.7551}&\underline{0.6642}& \underline{0.7221}&49&9 \\
		&\checkmark &\ding{53}&0.6710&0.7222&0.6206&0.6952 &13&4  \\
		&\ding{53}&\ding{53}&0.4895& 0.6176& 0.5279& 0.6367 &649 &9 \\
		\midrule
		\multirow{4}{*}{\makecell{MSRCV1}} &\checkmark &\checkmark&\bfseries0.7762&\bfseries0.7504&\bfseries0.6684&\bfseries0.7305&42 &42  \\
		&\ding{53}&\checkmark&0.6000& 0.6004&\underline{0.4907}& \underline{0.6145}&1045&42 \\
		&\checkmark &\ding{53}&\underline{0.6810}&\underline{0.6612}&0.4904&0.6132 &1060&35  \\
		&\ding{53}&\ding{53}&0.5619& 0.5493&0.3496& 0.5411 &2428 &42 \\
		\midrule
		\multirow{4}{*}{\makecell{Multi-omics}} &\checkmark &\checkmark&\bfseries0.9504&\bfseries0.7676&\bfseries0.8567&\bfseries0.9350&9 &9  \\
		&\ding{53}&\checkmark&\underline{0.9417}& \underline{0.7287}&\underline{0.8314}&\underline{0.9274}&1036&9 \\
		&\checkmark &\ding{53}&0.9359& 0.7115&0.8149&0.9154 &9&2  \\
		&\ding{53}&\ding{53}&0.9300& 0.6924 &0.7986& 0.9098&1813 &9  \\
		\bottomrule
	\end{tabular*}
\end{table*}

\subsubsection{Evaluation Indicators}
To assess the performance of clustering, four indicators are chosen as below. 
\begin{itemize} 
	\item Accuracy (ACC) \cite{zhan2018multiview}  measures the proportion of correctly clustered data points by comparing their assigned labels to the ground truth, which is defined by
	\begin{eqnarray}
	\operatorname{ACC}=\frac{\sum_{i=1}^n \mathbb{I}(r_i= \wp(c_i))}{n},
	\end{eqnarray}
	where  $r_i$ represents the ground-truth label of the $i$-th sample, $c_i$ denotes the corresponding learned clustering  assignment, and  $\wp(r_i)$ is the optimal permutation function that aligns predicted labels $c$ with the true labels $r$.
	\item Normalized mutual information (NMI) \cite{fang2023joint}  quantifies the statistical dependence between $r$ and $c$  similarity, which is given by
\begin{eqnarray}
	\operatorname{NMI}=\frac{\operatorname{MI}(r, c)}{\max\big\{\operatorname{EN}(r), \operatorname{EN}(c)\big\}},
\end{eqnarray}
where $\operatorname{MI}(r, c)$ and and $\operatorname{EN}(\cdot)$ return the mutual information and information entropy between $r$ and $c$, respectively.
\item Adjusted rand index (ARI) \cite{Wang2021} evaluates the agreement between $r$ and $c$ while accounting for the probability of random chance, which is calculated by 
\begin{eqnarray}
	\begin{aligned}
	\text{ARI}=
	\frac{\sum_{i, j=1}^K C_{n_{i, j}}^2-\frac{\sum_{i=1}^K C_{n_i^r}^2 \sum_{i=1}^K C_{n_i^c}^2}{C_n^2}}{\frac{1}{2}\left(\sum_{i=1}^K C_{n_i^r}^2+\sum_{i=1}^K C_{n_i^c}^2\right)-\frac{\sum_{i=1}^K C_{n_i^r}^2 \sum_{i=1}^K C_{n_i^c}^2}{C_n^2}},
	\end{aligned}
\end{eqnarray}
	where $K$ denotes the number of clusters,  $C_a^b$ calculates the count of possible combinations when selecting $b$ items from a total of $a$ items,  $n_{i,j}$ counts samples shared between $r$ and $c$ .
\item Fowlkes-Mallows index (FMI) \cite{Chen2020}  quantifies  clustering agreement by computing the geometric mean of pairwise precision and recall for sample pairs.
\begin{eqnarray}
	\mathrm{FMI}=\sqrt{\text { Precision } \times \text {Recall} },
\end{eqnarray}
where Precision reflects the accuracy in identifying negative instances, while Recall measures the completeness in detecting positive instances.
\end{itemize}

Generally, higher values of ACC, NMI, ARI, and FMI signify better clustering performance, reaching a maximum of 1 at perfect agreement where all samples are correctly clustered.

\subsection{Experimental Results}\label{4.2}

For better fitting, we leverage  the following distance metrics: the Euclidean distance for Gaussian-distributed data, the Manhattan distance for count data, as well as the Bernoulli log-likelihood for proportion-valued data or binary data. These metrics are empirically superior to alternative loss functions while maintaining computational efficiency. Simultaneously, we employ the adaptive parameter adjustment method described in \cite{Wang2021} to adjust parameters $\zeta_k$, $\eta$, $w_\iota$, $\beta$, and $\bar\beta^k$, and adopt the grid search method to select $\theta$ from $\{2^{-8}, 2^{-7}, 2^{-6},\ldots, 2^{0}\}$ and $\sigma$ from $\{0.1, 0.2, 0.3,\ldots, 1.5\}$. To be fair, all fusion-regularized clustering methods employ ADMM for optimization and utilize identical parameter selection protocols.

According to the characteristics of the real data shown in Figure \ref{fig4}, we select appropriate loss functions. 
As visualized in Table \ref{tab4},   our method has the highest clustering accuracy.  Boldface denotes the best  performance, while  underline indicates the second-best. In comparison with the three multi-view clustering methods, namely SCMvFS, CDMvFS,   and WMVEC, the average ARI of our method is enhanced by 11.23\%,  11.13\%, and 11.59\%, respectively. When compared with the three fusion regularized clustering methods (SCC, SGLCC, and iGecco+), our method increases the average ARI by approximately 16.12\%, 14.63\%, and 7.12\%, respectively. Among them,
the reason for the poor clustering effect of SCC and SGLCC is that they do not account for heterogeneity across individual views.  They uniformly treat all data as single-view data following a standard normal distribution.

To be more intuitive, Figure \ref{fig2} separately visualizes the clustering results of eight methods applied to  Dermatology, UCI-digits, and MSRCV1. Each subplot shows the t-SNE projections with color-coded cluster labels, demonstrating that our method achieves clearer separation.

\begin{figure*}[t]
	\centering
	\captionsetup[subfigure]{font=scriptsize, textfont=scriptsize, labelfont=scriptsize}
	\hspace*{-0.8cm} 
	\subfloat[(a) Authors]{
		\centering
		\includegraphics[width=3.1in]{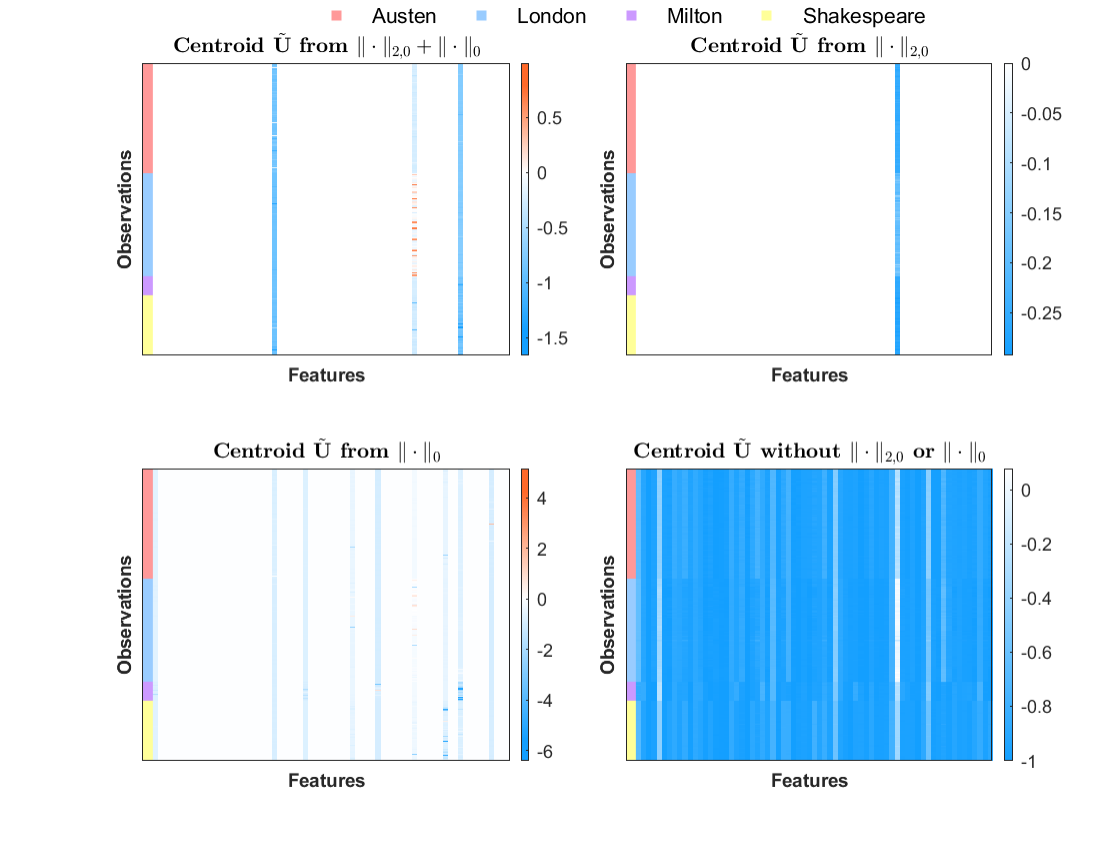}
		\label{fig:subfig5}
	}
	\subfloat[(b) Dermatology]{
		\centering
		\includegraphics[width=3.1in]{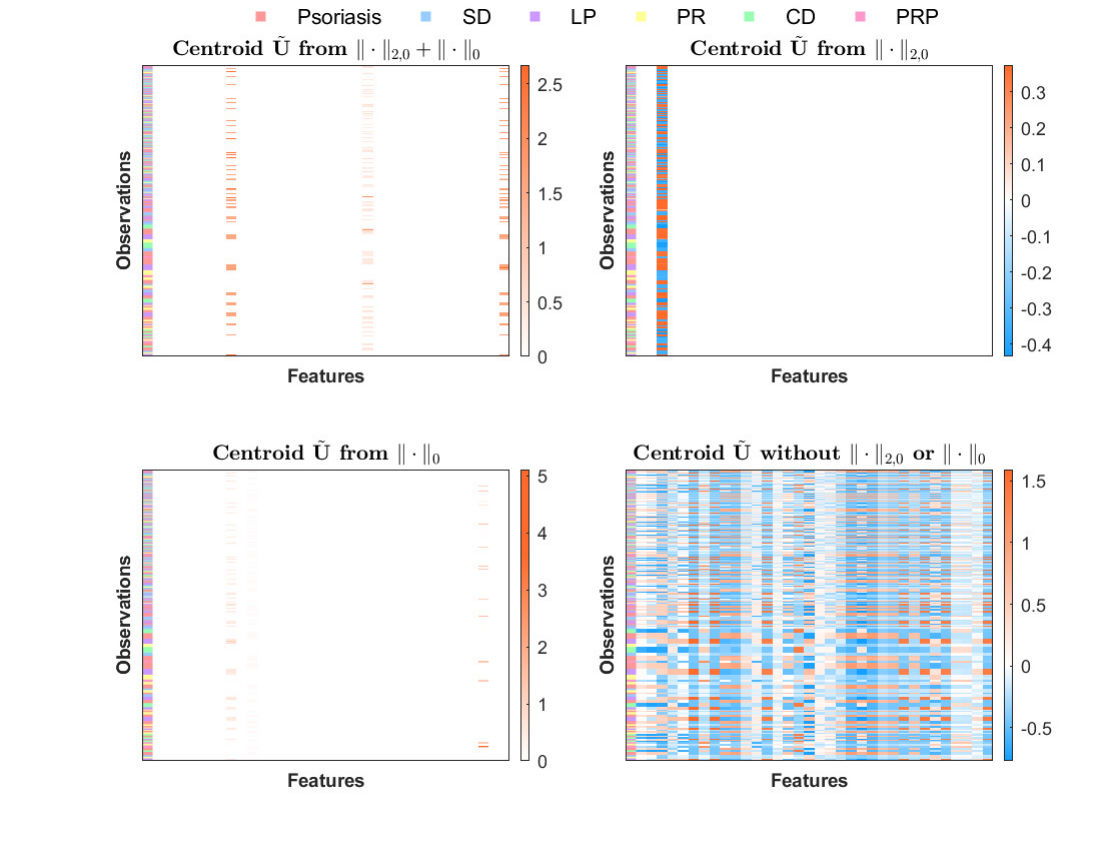}
		\label{fig:subfig6}
	}
	\caption{ Visualization of features selected through ablation experiments on  Authors and Dermatology. For each subfigure of each data,  the first column shows the ground-truth class distribution, and the remaining columns represent the selected features and their classification behavior under distinct feature selection mechanisms.}
	\label{fig7}
\end{figure*}

\subsection{Ablation Study}\label{4.3}

To further illustrate the effectiveness of the group sparsity strategy, we carry out sufficient ablation experiments for Authors, Lung-discrete, Dermatology,  UCI-digits, MSRCV1, and Multi-omics data.  
Let $T_{cur}$ and $T_{our}$ denote the feature sets selected by the current clustering method and  our clustering method, respectively. $|T_{cur}\cap T_{our}|$ represents the cardinality of the intersection of two sets $T_{cur}$ and $T_{our}$.
The intersection size  $|T_{cur}\cap T_{our}|$ quantifies their agreement.
Table \ref{tab6} presents a comparative evaluation of clustering performance, selected minimum feature number for optimal performance, and feature overlaps between our method and its three variants.

The quantitative comparison in Table \ref{tab6}, assessed through ACC, NMI, ARI and FMI metrics, reveals that our clustering algorithm consistently outperforms the three competing variants. Furthermore,
from Table \ref{tab6},   features selected by only considering $\|\cdot\|_{2,0}$ or $\|\cdot\|_{0}$ is inconsistent with that by our method. In particular, for Lung-discrete and MSRCV1, the feature sets selected by these two methods both satisfy $|T_{cur}\cap T_{our}|<|T_{our}|$. When neither module is added, all features will be considered important.
By analyzing both selected feature quantities and intersections from the aforementioned experiments, we can find that  $\|\cdot\|_{2,0}$ and $\|\cdot\|_{0}$ induce divergent feature selection behaviors, yielding  fundamentally different impacts on model performance.  Visualization of the ablation experiments on the Authors and Dermatology data clearly illustrates this phenomenon. As depicted  in Figure \ref{fig7}, 
even when group sparsity ($\|\cdot\|_{2,0}+\|\cdot\|_{0}$) and intra-group sparsity ($\|\cdot\|_{0}$) select the same number of features on Dermatology, they share only one common feature. Notably, features selected by group/intra-group sparsity exhibit partial activation patterns, whereas those from inter-group sparsity show full activation. The variant without $\|\cdot\|_{2,0}$ or $\|\cdot\|_{0}$ fails to perform feature selection.
In summary, different components $\|\cdot\|_{2,0}$ and  $\|\cdot\|_{0}$ play certain roles in feature selection and clustering accuracy, and our model with group sparsity has a distinct advantage.

\subsection{Discussions} \label{4.4}

\subsubsection{Stability Analysis}

\begin{figure*}[t]
	\centering
	\captionsetup[subfigure]{font=scriptsize, textfont=scriptsize, labelfont=scriptsize}
	\hspace*{-0.3cm} 
	\subfloat[(a1) Authors (ACC)]{
		\centering
		\includegraphics[width=0.24\textwidth]{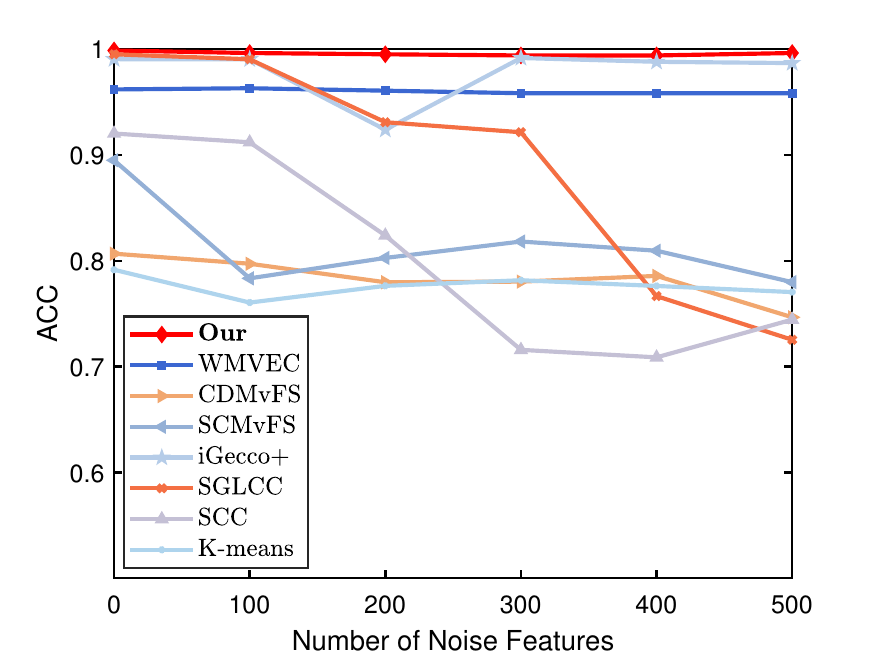}
		\label{fig:subfig01}
	}	\hspace*{-0.3cm} 
	\subfloat[(a2) Authors (NMI)]{
		\centering
		\includegraphics[width=0.24\textwidth]{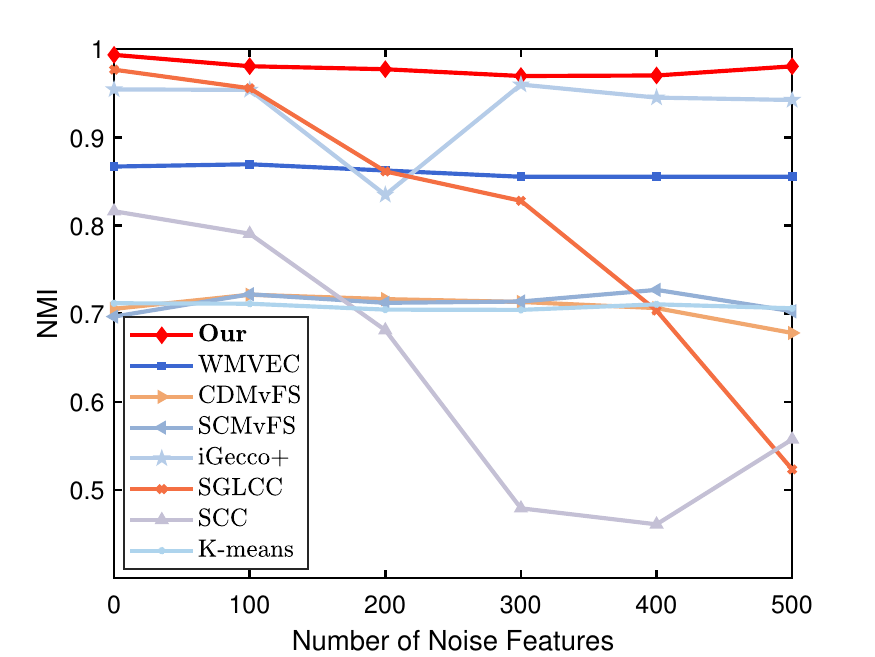}
		\label{fig:subfig02}
	}	\hspace*{-0.3cm} 
	\subfloat[(a3) Authors (ARI)]{
		\centering
		\includegraphics[width=0.24\textwidth]{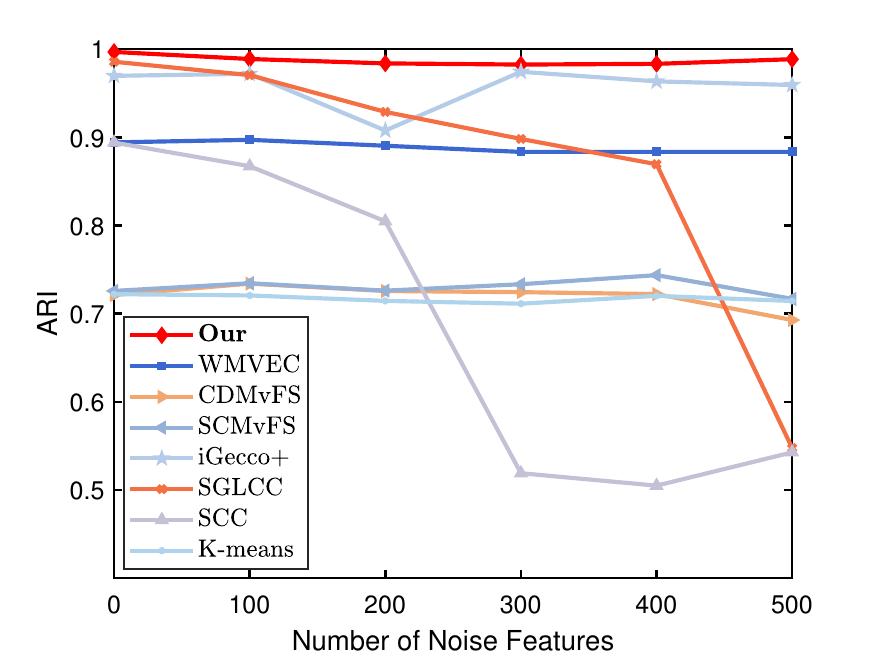}
		\label{fig:subfig03}
	}	\hspace*{-0.3cm} 
	\subfloat[(a4) Authors (FMI)]{
		\centering
		\includegraphics[width=0.24\textwidth]{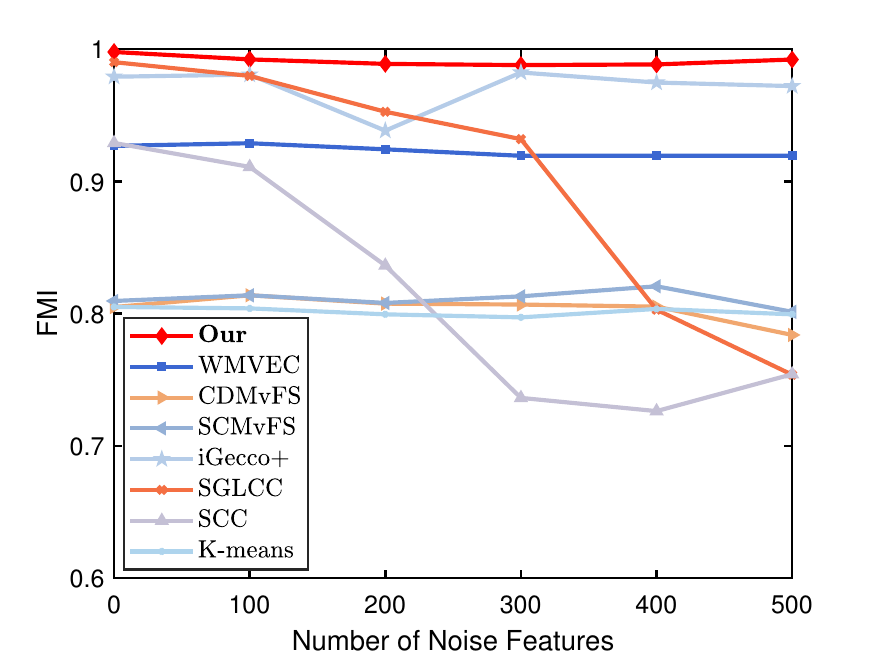}
		\label{fig:subfig04}
	}
	\\
	\hspace*{-0.3cm} 
	\subfloat[(b1) Dermatology (ACC)]{
		\centering
		\includegraphics[width=0.24\textwidth]{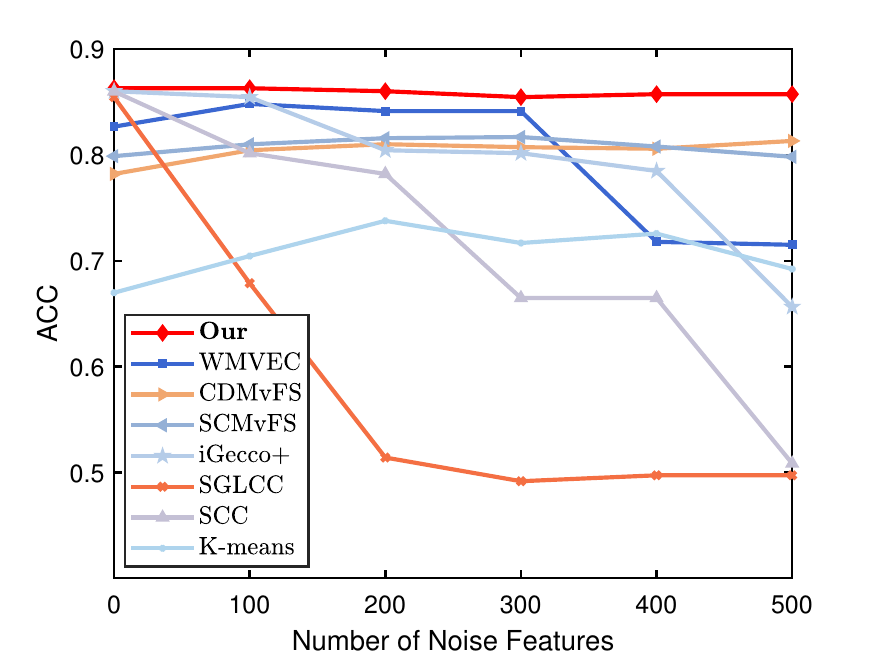}
		\label{fig:subfig05}
	}	\hspace*{-0.3cm} 
	\subfloat[(b2) Dermatology (NMI)]{
		\centering
		\includegraphics[width=0.24\textwidth]{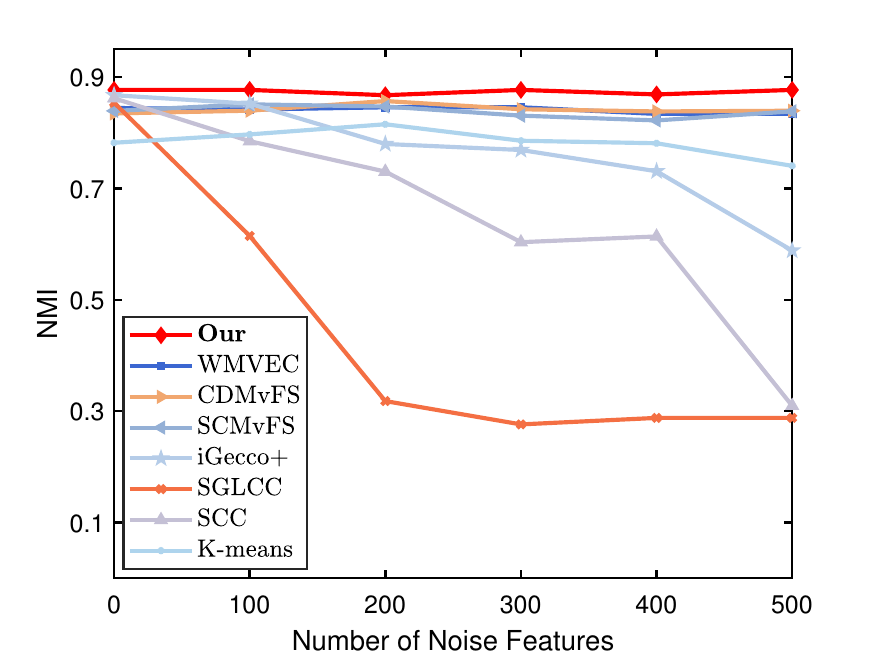}
		\label{fig:subfig06}
	}	\hspace*{-0.3cm} 
	\subfloat[(b3) Dermatology (ARI)]{
		\centering
		\includegraphics[width=0.24\textwidth]{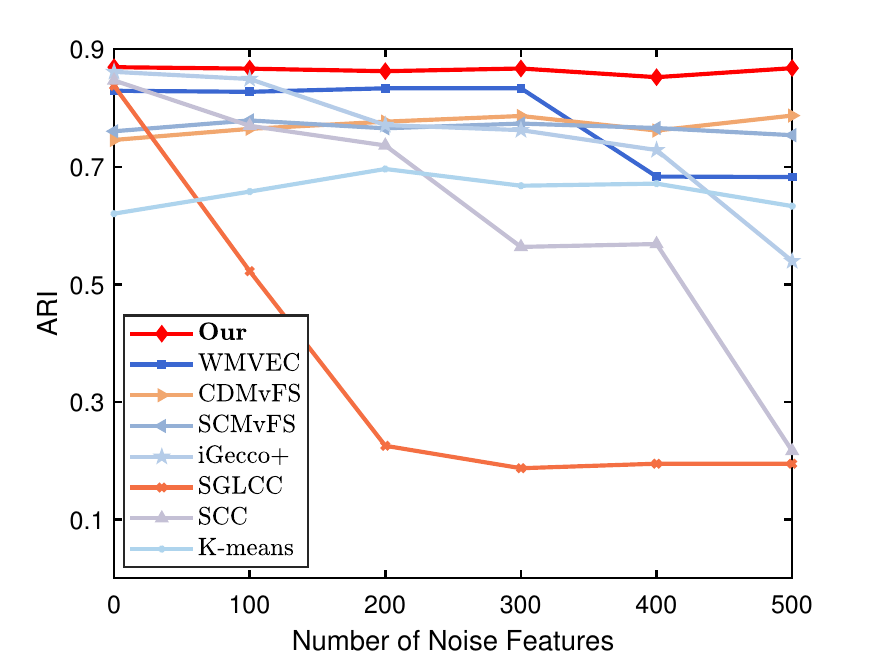}
		\label{fig:subfig07}
	}	\hspace*{-0.3cm} 
	\subfloat[(b4) Dermatology (FMI)]{
		\centering
		\includegraphics[width=0.24\textwidth]{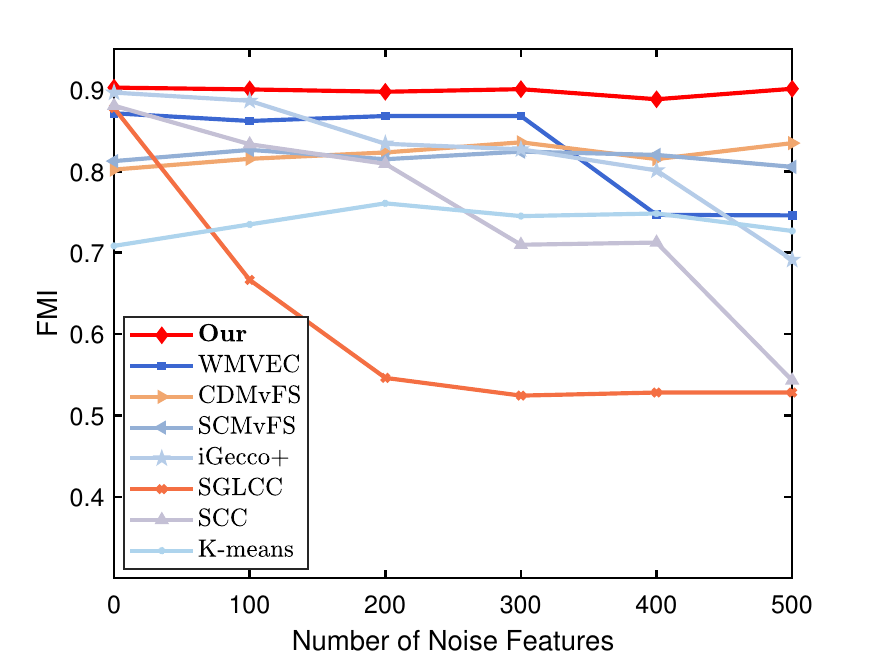}
		\label{fig:subfig08}
	}	
	\caption{Visual comparison of clustering results from 8 methods on Authors and Dermatology as noise levels increase.}
	\label{fig8}
\end{figure*}

In this study, we corrupt the non-Gaussian Authors data with Poisson-distributed noise. For Dermatology data, we inject equal amounts of Poisson-distributed noise for the first view (clinical assessments) and Bernoulli-distributed noise for the other view (histopathological examinations). A line chart is employed to compare the clustering performance between our approach and eight competing methods. As evident from Figure \ref{fig8}, our method consistently demonstrates stable and superior performance as the number of noisy features increases. These results further suggest the robustness  and effectiveness of our proposed approach.
\begin{figure*}[t]
	\centering
	\captionsetup[subfigure]{font=scriptsize, textfont=scriptsize, labelfont=scriptsize}
	\hspace*{-0.8cm} 
	\subfloat[(a) Authors]{
		\centering
		\includegraphics[width=3.53in]{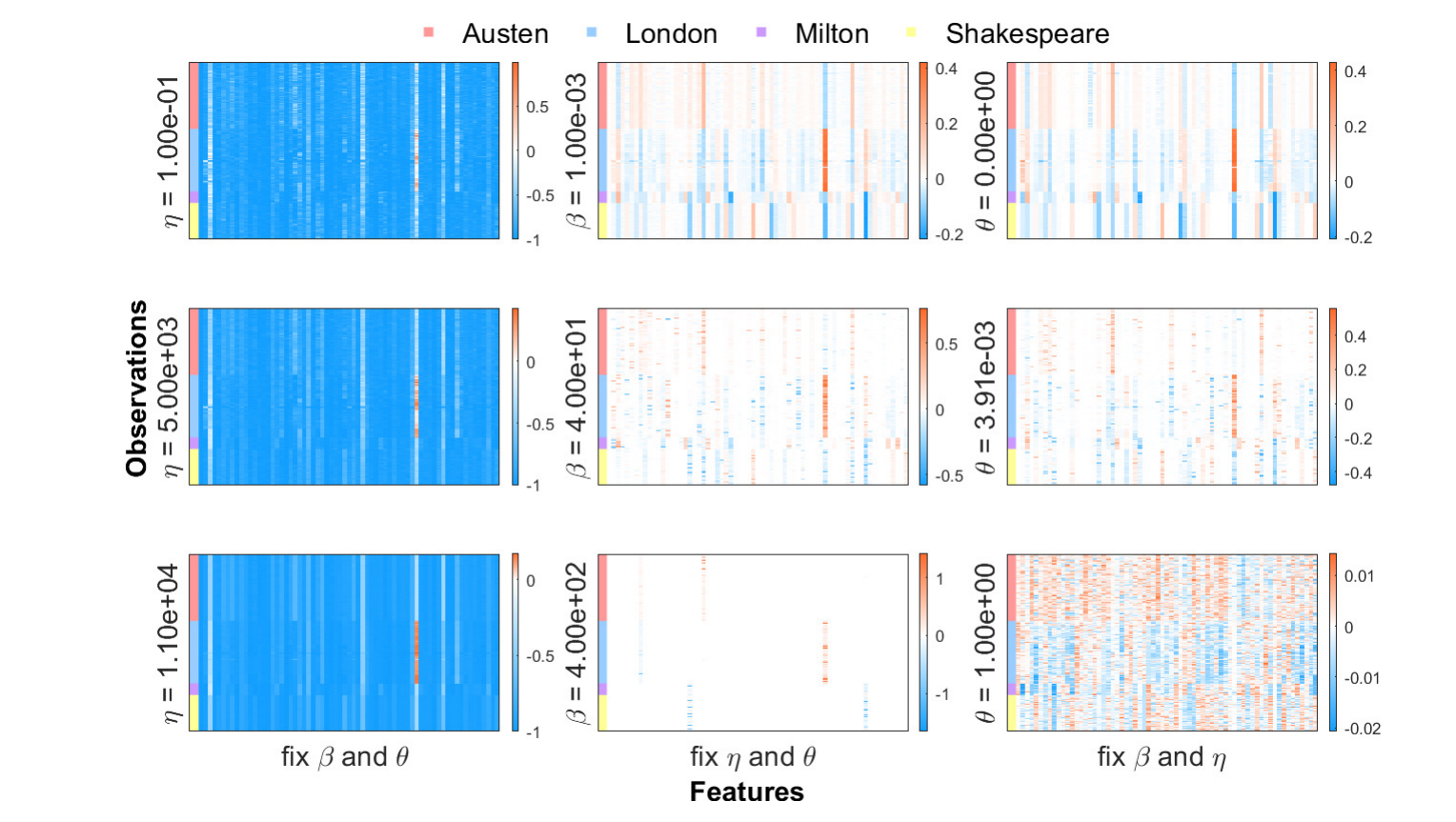}
		\label{fig:subfig1}
	}\hspace*{-8mm}
	\subfloat[(b) Dermatology]{
		\centering
		\includegraphics[width=3.53in]{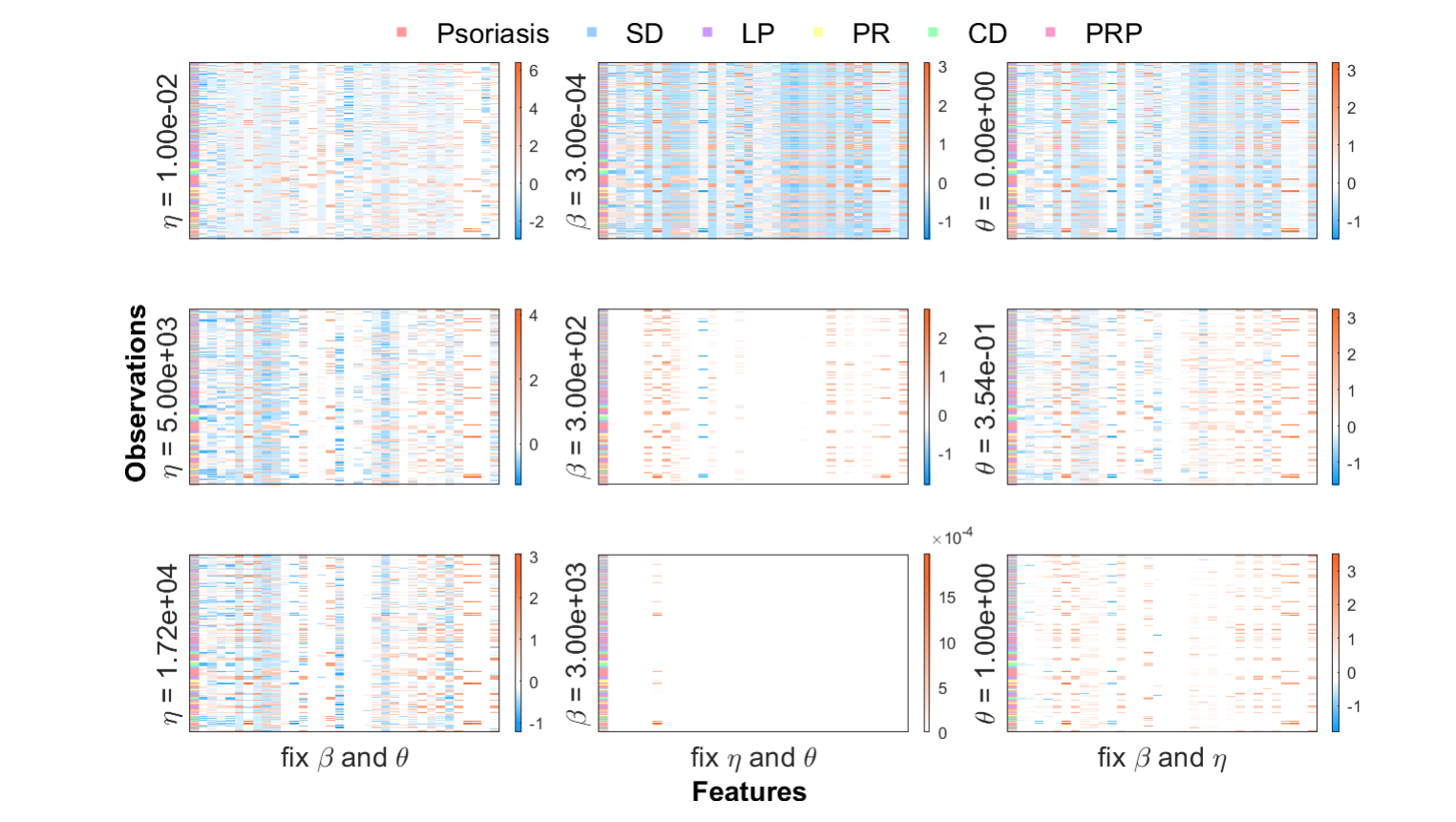}
		\label{fig:subfig2}
	}
	\caption{Clustering path $\eta$ and regularization path $\beta$, $\theta$ of our method  on Authors and Dermatology. For the subplots of each data, from left to right, we characterize the responses of $\bm{\tilde{U}}$  by fixing $\beta$ and $\theta$ and adjusting $\eta$, the responses of $\bm{\tilde{U}}- \bm{\bar{X}}$  by fixing $\eta$ and $\theta$ and adjusting $\beta$, and the responses of $\bm{\tilde{U}}- \bm{\bar{X}}$  by fixing $\beta$ and $\eta$ and adjusting $\theta$, respectively.}
	\label{fig.4}
\end{figure*}

\subsubsection{Parameter Analysis}

We perform sensitivity analysis of the tuning parameters $\eta$, $\beta$, and $\theta$ by applying  our method on Authors and Dermatology. Figure \ref{fig.4} reflects how $\bm{\tilde{U}}$ and  $\bm{\tilde{U}}-\bm{\bar{X}}$ change as a function of three tuning parameters, where the size of  $\bm{\tilde{U}}$ and  $\bm{\tilde{U}}- \bm{\bar{X}}$ can clearly describe the situation of clustering and feature selection  respectively. It is found that
\begin{itemize}
	\item If we fix  $\beta$ and $\theta$,  the number of clusters decreases as the increase of $\eta$. 
	\item If we fix  $\eta$ and $\theta$,  the number of selected features decreases with the increase of $\beta$. 
	\item  Fixing $\beta$ and $\eta$,  when $\theta=0$, whole  informative features are selected. As $\theta$ increases, local informative features are selected. Until $\theta=1$, the model does not have feature selection capability because all features are selected. 
\end{itemize}

\begin{figure*}[t]
	\centering
	\captionsetup[subfigure]{font=scriptsize, textfont=scriptsize, labelfont=scriptsize}
	\hspace*{-0.7cm}
	\subfloat[(a) Authors]{
		\centering
		\includegraphics[width=3.4in]{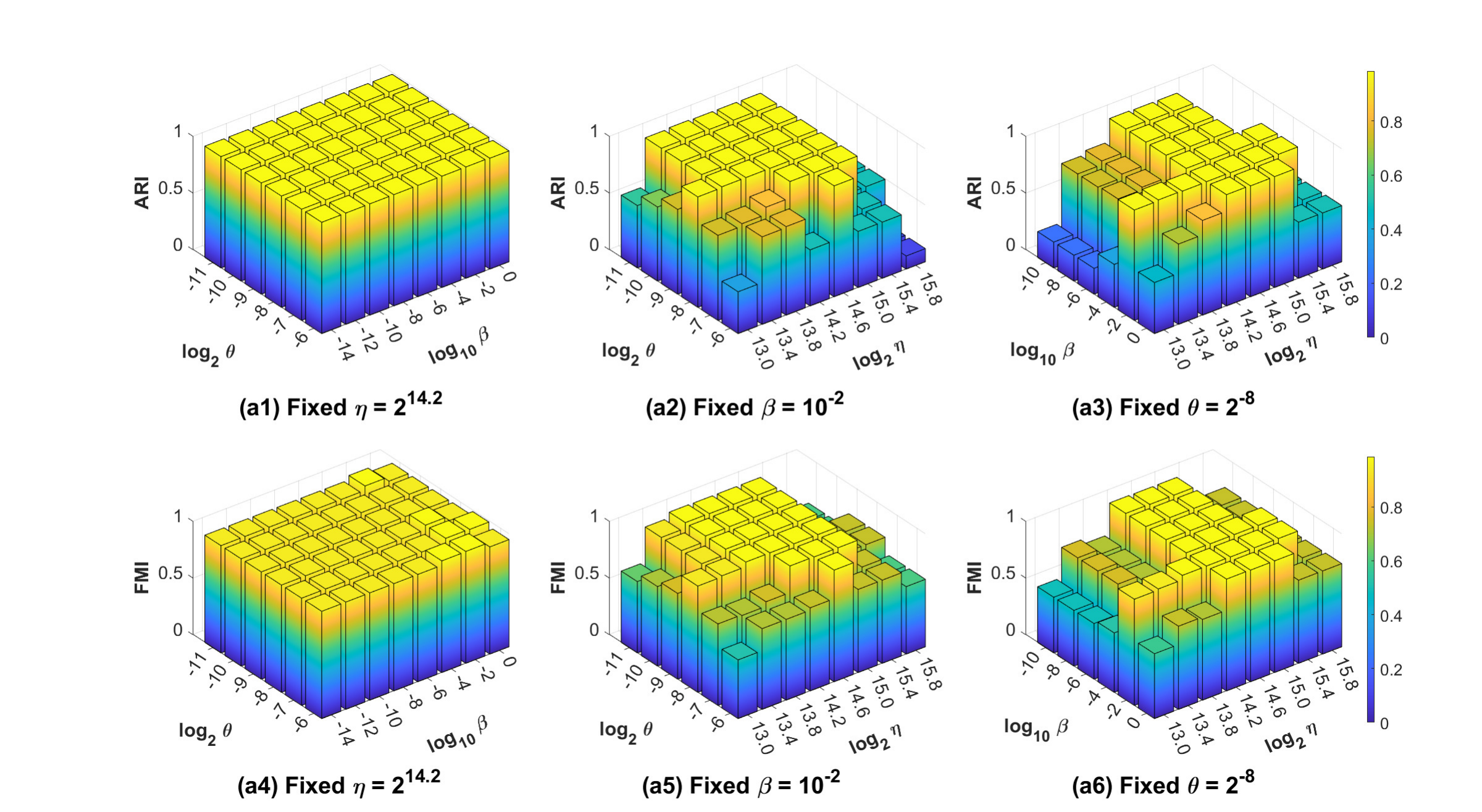}
		\label{fig:subfig3}
	}\hspace*{-7.9mm}
	\subfloat[(b) Dermatology]{
		\centering
		\includegraphics[width=3.4in]{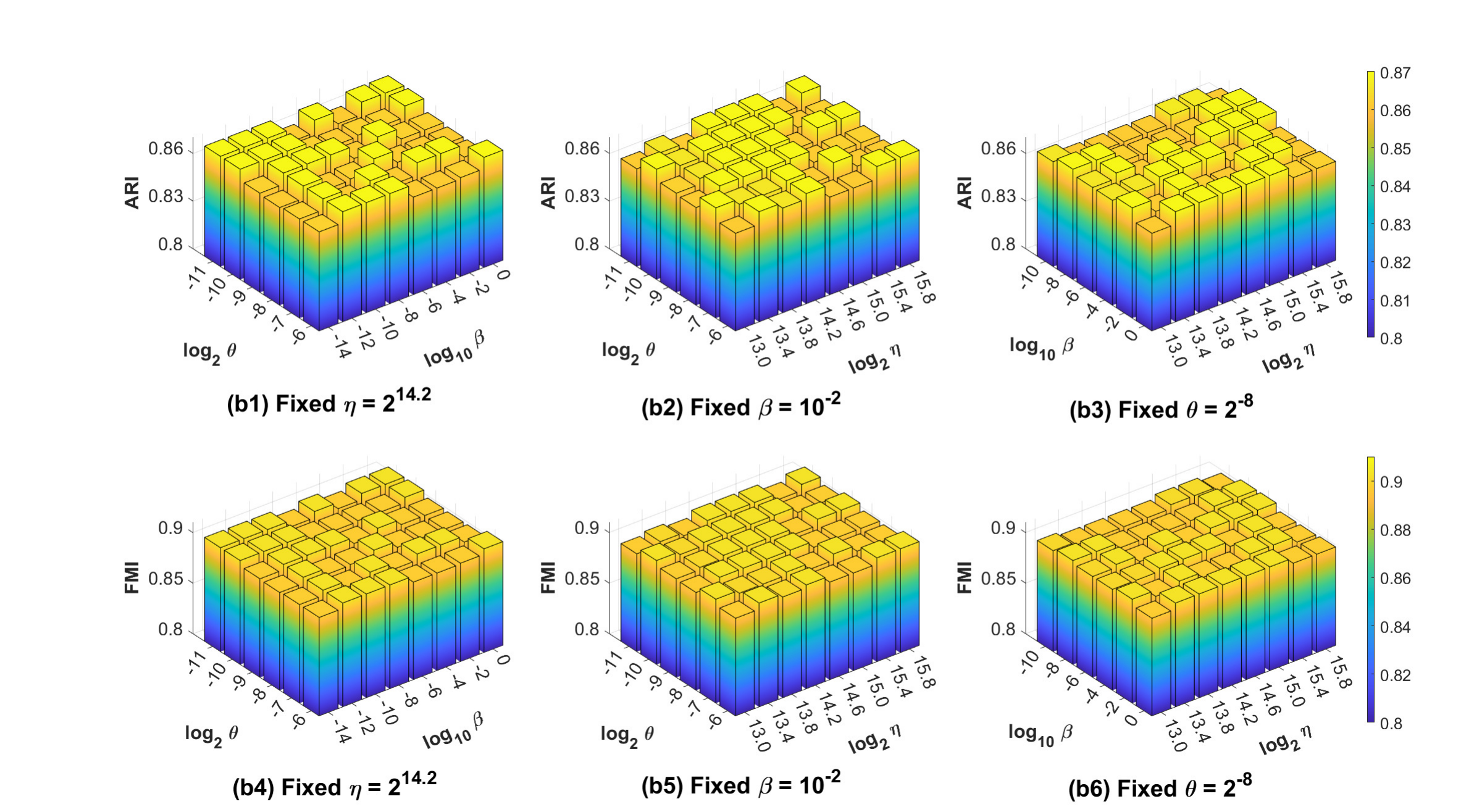}
		\label{fig:subfig4}
	}
	\caption{Parameter sensitivity using ARI and FMI on Authors and Dermatology.}
	\label{fig5}
\end{figure*}

In addition, Figure \ref{fig5} shows the 3D histograms of the change in ARI and FMI.  It can be observed that, for Authors data, when $\eta$ is fixed, the ARI and FMI are not sensitive to changes in the penalty parameters $\theta$ and $\beta$; When $\beta$  or $\theta$ is fixed,  our proposed method maintains good performance over a certain range.  In contrast, the Dermatology data exhibits extremely low sensitivity to all parameter changes, with virtually negligible effects. 
\begin{figure*}[t]
	\centering
	\captionsetup[subfigure]{font=scriptsize, textfont=scriptsize, labelfont=scriptsize}
	\subfloat[(a) Authors]{
		\centering
		\includegraphics[width=1.5in]{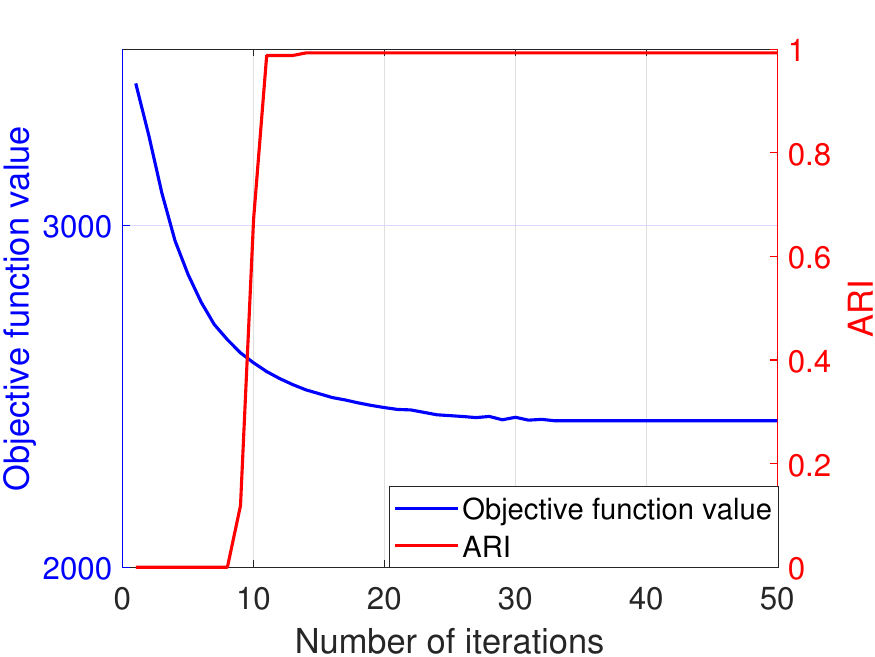}
		\label{fig:subfig11}
	}\hspace*{3mm}
	\subfloat[(b) Lung-discrete]{
		\centering
		\includegraphics[width=1.5in]{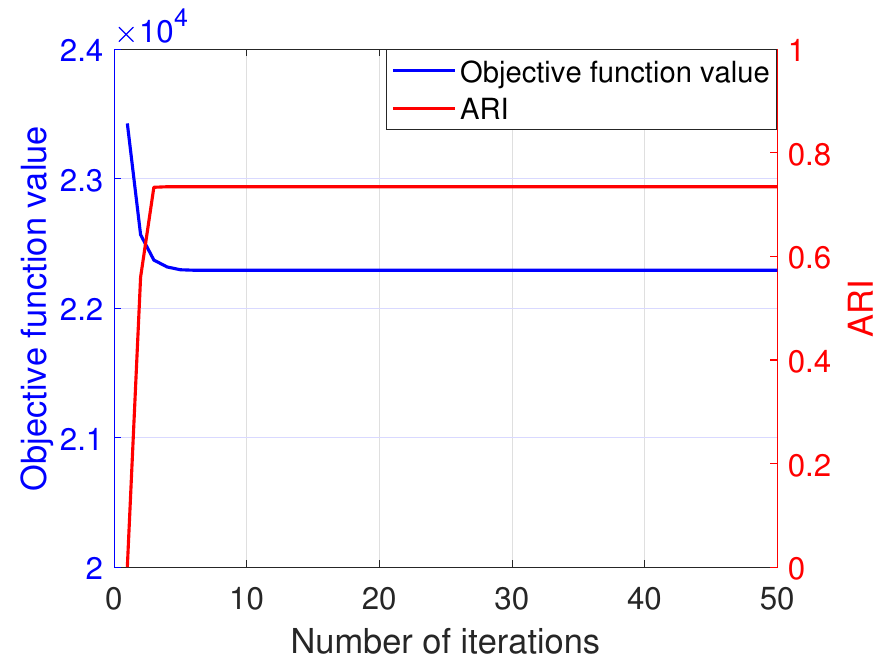}
		\label{fig:subfig12}
	}\hspace*{3mm}
	\subfloat[(c) Dermatology]{
		\centering
		\includegraphics[width=1.5in]{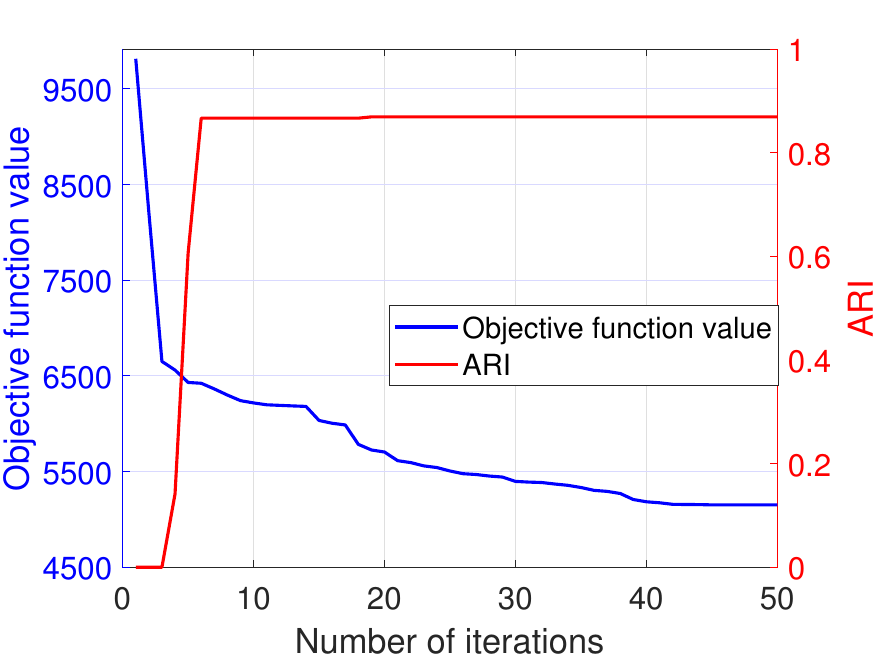}
		\label{fig:subfig13}
	}	\\
	\subfloat[(d) UCI-digits]{
		\centering
		\includegraphics[width=1.5in]{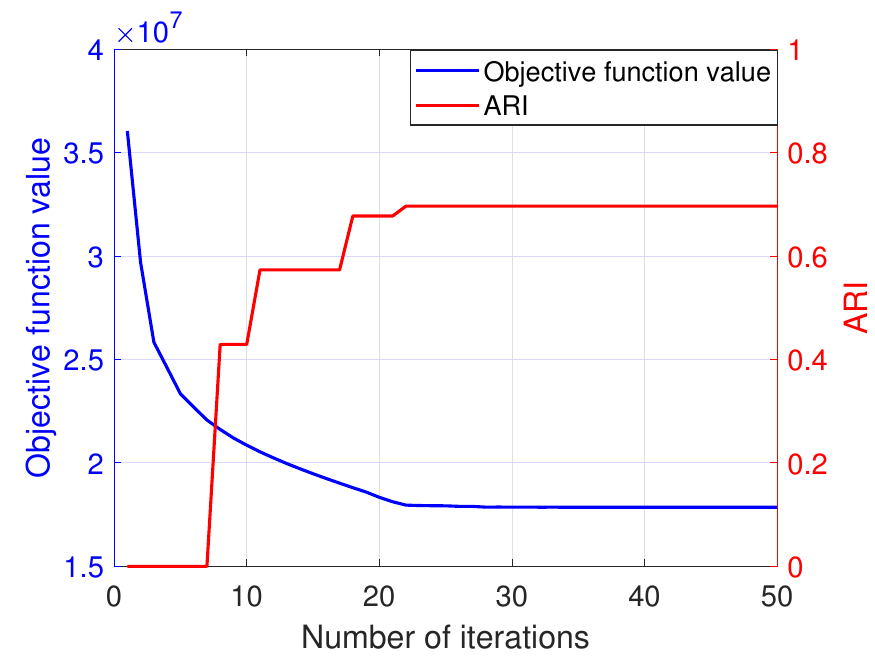}
		\label{fig:subfig14}
	}\hspace*{3mm}
	\subfloat[(e) MSRCV1]{
		\centering
		\includegraphics[width=1.5in]{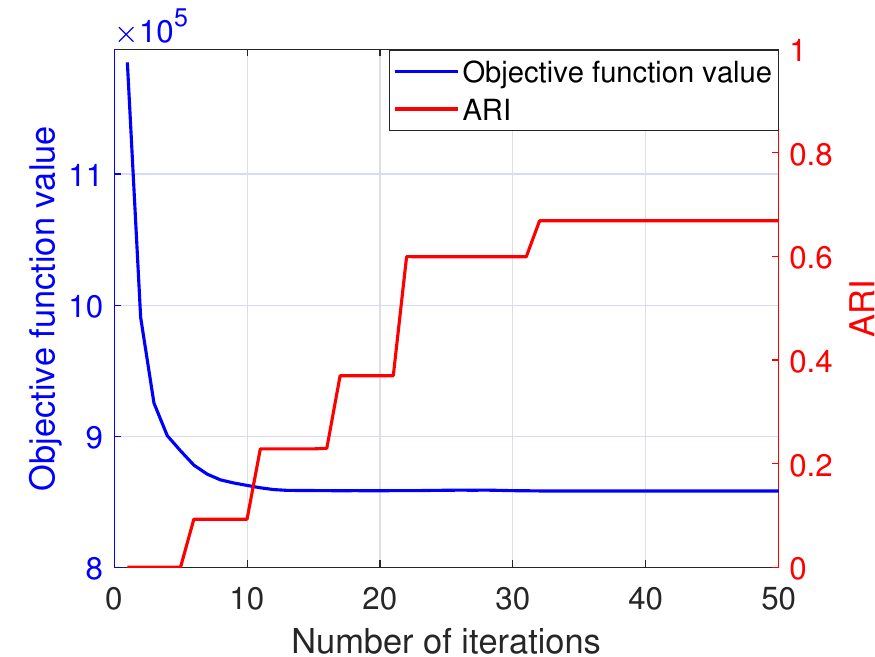}
		\label{fig:subfig15}
	}\hspace*{3mm}
	\subfloat[(f) Multi-omics]{
		\centering
		\includegraphics[width=1.5in]{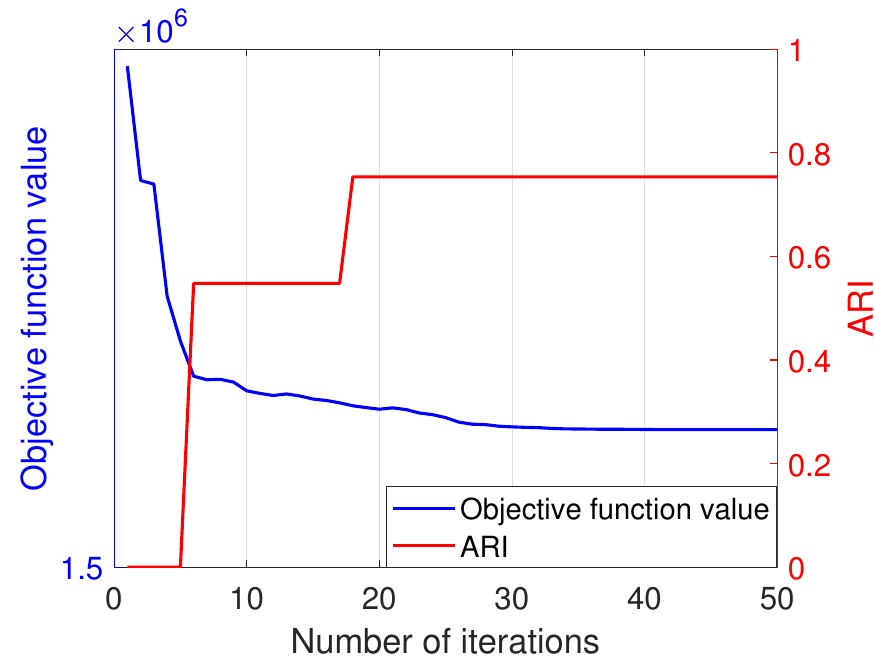}
		\label{fig:subfig17}
	}
	\caption{Objective function value and ARI versus  iteration number  of our method on  6 real data.}
	\label{fig6}
\end{figure*}

\subsubsection{Convergence}

Figure \ref{fig6} displays the evolution of both the objective function value in (\ref{XX01}) and ARI versus iteration numbers on six real data. It is evident from the plot that  the objective function value exhibits an overall decreasing trend with iterations, while ARI shows a corresponding increase, with both metrics ultimately converging to stable values.

\section{Conclusions}\label{section6}

In this paper, we propose a novel multi-view fusion regularized clustering method to address \emph{heterogeneity} and \emph{redundancy}.  To the best of our knowledge, this is the first time that non-convex group sparsity is introduced into fusion regularized clustering methods. In addition, we develop an efficient  ADMM algorithm.  Extensive experimental results verify its advantages, especially on the Multi-omics dataset, where the proposed method improves ACC by at least 1.75\%, NMI by at least 4.89\%, ARI by at least 5.01\%, and FMI by at least 2.01\% over the compared methods.

In future work,  we will focus on  developing high-performance algorithms to accelerate computation. Besides, combining clustering with deep learning to improve representation performance is also a potential  direction.

%
%
\printcredits

\bibliographystyle{model1-num-names}
\bibliography{mybib}

\end{document}